\documentclass[aos,citesort,dvips]{arximspdf}
\usepackage{subenv}
\usepackage{accents}
\usepackage{graphicx}

\doi{10.1214/09-AOS691}
\volume{38}
\issue{3}
\pubyear{2010}
\firstpage{1287}
\lastpage{1319}

\makeatletter

\newcommand{\twobar}[1]{\accentset{\vspace*{-1pt}=}{#1}}
\newcommand{\onebar}[1]{\accentset{\vspace*{-2pt}-}{#1}}

\newcommand{\nnorm}{|\!|\!|}

\newtheorem{theos}{Theorem}
\newtheorem{props}{Proposition}
\newtheorem{lems}{Lemma}
\newtheorem{cors}{Corollary}

\makeatother

\begin{document}
\begin{frontmatter}

\title{High-dimensional Ising model selection using ${\ell_1}$-regularized logistic regression}
\runtitle{High-dimensional Ising model selection}

\begin{aug}
\author[A]{\fnms{Pradeep} \snm{Ravikumar}\corref{}\thanksref{t1,t2,t3}\ead[label=e1]{pradeepr@stat.berkeley.edu}},
\author[A]{\fnms{Martin J.} \snm{Wainwright}\thanksref{t3}\ead[label=e2]{wainwrig@stat.berkeley.edu}}\\ and
\author[B]{\fnms{John D.} \snm{Lafferty}\thanksref{t1}\ead[label=e3]{lafferty@cs.cmu.edu}}
\runauthor{P. Ravikumar, M. J. Wainwright and J. D. Lafferty}
\affiliation{University of California, Berkeley, University of
California, Berkeley and Carnegie Mellon University}
\address[A]{P. Ravikumar\\
M. J. Wainwright\\
Department of Statistics\\
University of California\\
Berkeley, California 94720\\
USA\\
\printead{e1}\\
\phantom{E-mail: }\printead*{e2}} 
\address[B]{J. D. Lafferty\\
Computer Science Department\\
\quad and Machine Learning Department\\
Carnegie Mellon University\\
Pittsburgh, Pennsylvania 15213\\
USA\\
\printead{e3}}
\end{aug}

\thankstext{t1}{Supported in part by NSF Grants IIS-0427206 and
CCF-0625879.}

\thankstext{t2}{Supported in part by a Siebel Scholarship.}

\thankstext{t3}{Supported in part by NSF Grants DMS-06-05165 and CCF-0545862.}

\received{\smonth{10} \syear{2008}}
\revised{\smonth{1} \syear{2009}}

%
\begin{abstract}
We consider the problem of estimating the graph associated with a
binary Ising Markov random field. We describe a method based on
$\ell_1$-regularized logistic regression, in which the neighborhood of
any given node is estimated by performing logistic regression subject
to an $\ell_1$-constraint. The method is analyzed under
high-dimensional scaling in which both the number of nodes $p$ and
maximum neighborhood size $d$ are allowed to grow as a function of the
number of observations $n$. Our main results provide sufficient
conditions on the triple $(n, p, d)$ and the model parameters for the
method to succeed in consistently estimating the neighborhood of every
node in the graph simultaneously. With coherence conditions imposed on
the population Fisher information matrix, we prove that consistent
neighborhood selection can be obtained for sample sizes $n = \Omega(d^3
\log p)$ with exponentially decaying error.
When these same conditions are imposed directly on the sample matrices,
we show that a reduced sample size of $n = \Omega(d^2 \log p)$ suffices
for the method to estimate neighborhoods consistently. Although this
paper focuses on the binary graphical models, we indicate how a
generalization of the method of the paper would apply to general
discrete Markov random fields.
\end{abstract}

%

\end{frontmatter}

\section{Introduction}\label{sec1}

Undirected graphical models, also known as Markov random fields, are
used in a variety of domains, including statistical physics
\cite{Ising25}, natural language processing \cite{Manning99}, image
analysis \cite{Woods78,Hassner80,Cross83} and spatial statistics
\cite{Ripley81}, among others. A Markov random field (MRF) is specified
by an undirected graph $G= (V, E)$ with vertex set $V= \{1, 2, \ldots,
p\}$ and edge set $E\subset V\times V$. The structure of this graph
encodes certain conditional independence assumptions among subsets of
the $p$-dimensional discrete random variable $X= (X_1, X_2,\ldots,
X_p)$ where variable $X_i$ is associated with vertex $i \in V$. One
important problem for such models is to estimate the underlying graph
from $n$ independent and identically distributed samples $\{ x^{(1)},
x^{(2)}, \ldots, x^{(n)} \}$ drawn from the distribution specified by
some Markov random field. As a concrete illustration, for binary random
variables, each vector-valued sample $x^{(i)} \in\{0,1\}^p$ might
correspond to the votes of a set of $p$ politicians on a particular
bill, and estimating the graph structure amounts to detecting
statistical dependencies in these voting patterns (see Banerjee, Ghaoui
and d'Aspr\'{e}mont \cite{BanGhaAsp08} for further discussion of this
example).

Due to both its importance and difficulty, the problem of structure
learning for discrete graphical models has attracted considerable
attention. The absence of an edge in a graphical model encodes a
conditional independence assumption. Constraint-based approaches
\cite{spirtes00} estimate these conditional independencies from the
data using hypothesis testing and then determine a graph that most
closely represents those independencies. Each graph represents a model
class of graphical models; learning a graph then is a model class
selection problem. Score-based approaches combine a metric for the
complexity of the graph with a measure of the goodness of fit of the
graph to the data; for instance, log-likelihood of the maximum
likelihood parameters given the graph, to obtain a \textit{score} for
each graph. The score is used together with a search procedure that
generates candidate graph structures to be scored. The number of graph
structures grows super-exponentially, however, and Chickering
\cite{chickering95} shows that this problem is in general NP-hard.

A complication for undirected graphical models involving discrete
random variables is that typical score metrics involve the partition
function or cumulant function associated with the Markov random field.
For general undirected MRFs, calculation of this partition
function is computationally intractable \cite{Welsh93}. The space of
candidate structures in scoring based approaches is thus typically
restricted to either directed graphical models \cite{dasgupta99} or
to simple sub-classes of undirected graphical models such as those
based on trees \cite{chowliu68} and hypertrees \cite{srebro03}.
Abbeel, Koller and Ng \cite{AbbKolNg06} propose a method for learning factor
graphs based on local conditional entropies and thresholding and
analyze its behavior in terms of Kullback--Leibler divergence between
the fitted and true models. They obtain a sample complexity that
grows logarithmically in the number of vertices $p$, but the computational
complexity grows at least as quickly as ${\mathcal{O}}(p^{d+1})$
where $d$ is the maximum neighborhood size in the graphical
model. This order of complexity arises from the fact that for each
node, there are ${p\choose d} = {\mathcal{O}}(p^d)$
possible neighborhoods of size $d$ for a graph with $p$
vertices. Csisz\'{a}r and Talata \cite{Csiszar06} show consistency of a
method that uses pseudo-likelihood and a modification of the BIC
criterion, but this also involves a prohibitively expensive search.

The main contribution of this paper is a careful analysis of the
computational and statistical efficiency of a simple method for
graphical model selection. The basic approach is straightforward: it
involves performing $\ell_1$-regularized logistic regression of each
variable on the remaining variables, and then using the sparsity
pattern of the regression vector to infer the underlying neighborhood
structure. Our analysis is high-dimensional in nature, meaning that
both the model dimension $p$ as well as the maximum neighborhood
size $d$ may tend to infinity as a function of the size
$n$. Our main result shows that under mild assumptions on
the population Fisher information matrix, consistent neighborhood
selection is possible using $n= \Omega(d^3 \log p)$
samples and computational complexity
${\mathcal{O}}(\max\{n, p\} p^3 )$. We also show that when
the same assumptions are imposed
directly on the sample matrices, $n = \Omega(d^2 \log p)$
samples suffice for consistent neighborhood selection with the same
computational complexity. We focus in this paper on
binary Ising models, but indicate in Section \ref{SecGenMRF} a
generalization of the method applicable to general discrete Markov
random fields.

The technique of $\ell_1$-regularization for estimation of sparse
models or signals has a long history in many fields (for instance,
see \cite{Tropp06} for one survey). A surge of recent work has shown
that $\ell_1$-regularization can lead to practical algorithms with
strong theoretical guarantees
(e.g., \cite
{CandesTao06,DonohoElad03,Meinshausen06,ng04,Tropp06,Wainwright06new,ZhaoYu06}).
Despite the well-known computational intractability of computing
marginals and likelihoods for discrete MRFs \cite{Welsh93}, our method
is computationally efficient; it involves neither computing the
normalization constant (or partition function) associated with the
Markov random field nor a combinatorial search through the space of
graph structures. Rather, it requires only the solution of standard
convex programs with an overall computational complexity of order
${\mathcal{O}}(\max\{p, n\} p^3)$ and is thus well suited
to high-dimensional problems \cite{KohKimBoy07}. Conceptually, like
the work of Meinshausen and B\"{u}hlmann \cite{Meinshausen06} on
covariance selection in Gaussian graphical models, our approach can be
understood as using a type of pseudo-likelihood based on the local
conditional likelihood at each node. In contrast to the Gaussian
case, where the exact maximum likelihood estimate can be computed
exactly in polynomial time, this use of a surrogate loss function is
essential for discrete Markov random fields given the intractability
of computing the exact likelihood \cite{Welsh93}.

Portions of this work were initially reported in a conference
publication \cite{WaiRavLaf06}, with the weaker result that
$n = \Omega(d^6 \log d+ d^5 \log p)$ samples
suffice for consistent Ising model selection.
Since the appearance of that paper, other researchers have also studied
the problem of model selection in
discrete Markov random fields. For the special case of bounded degree
models, Bresler, Mossel and Sly \cite{Bresler08} describe a simple search-based
method, and prove under relatively mild assumptions that it can
recover the graph structure with $\Theta(\log p)$ samples.
However, in the absence of additional restrictions, the computational
complexity of the method is ${\mathcal{O}}(p^{d+1})$. In other
work, Santhanam and Wainwright \cite{SanWai08} analyze the
information-theoretic limits of graphical model selection, providing
both upper and lower bounds on various model selection procedures, but
these methods also have prohibitive computational costs.

The remainder of this paper is organized as follows. We begin in
Section \ref{SecBackground} with background on discrete graphical
models, the model selection problem and logistic regression. In Section
\ref{SecOutline}, we state our main result, develop some of its
consequences and provide a high-level outline of the proof. Section
\ref{SecFixedDesign} is devoted to proving a result under stronger
assumptions on the sample Fisher information matrix whereas Section
\ref{SecPopulation} provides concentration results linking the
population matrices to the sample versions. In Section
\ref{SecExperiments}, we provide some experimental results that
illustrate the practical performance of our method and the close
agreement between theory and practice. Section \ref{SecGenMRF}
discusses an extension to more general Markov random fields, and we
conclude in Section \ref{SecDiscussion}.

\subsection*{Notation} For the convenience of the reader, we
summarize here notation to be used throughout the paper. We use the
following standard notation for asymptotics: we write $f(n) =
{\mathcal{O}}(g(n))$ if $f(n) \leq Kg(n)$ for some constant $K<
\infty$, and $f(n) = \Omega(g(n))$ if $f(n) \geq K' g(n)$ for
some constant $K' > 0$. The notation $f(n) = \Theta(g(n))$
means that $f(n) = {\mathcal{O}}(g(n))$ and $f(n) = \Omega(g(n))$.
Given a
vector $v \in{\mathbb{R}}^d$ and parameter $q \in[1, \infty]$, we use
$\|v\|_q$ to denote the usual $\ell_q$ norm. Given a matrix $A
\in{\mathbb{R}}^{a \times b}$ and parameter $q \in[1, \infty]$, we use
$\nnorm A \nnorm_{{q}}$ to denote the induced matrix-operator norm
with $A$ viewed as a mapping from $\ell_q^b \rightarrow\ell
_q^a$ (see Horn
and Johnson \cite{Horn85}). Two examples of particular importance in
this paper are the spectral norm $\nnorm A \nnorm_{{2}}$,
corresponding to the maximal singular value of $A$, and the
$\ell_\infty$ matrix norm, given by $\nnorm A \nnorm_{{\infty}}
=\max
_{j=1, \ldots, a} {\sum_{k=1}^b} |A_{jk}|$. We make use
of the bound $\nnorm A \nnorm_{{\infty}} \leq\sqrt{a}
\nnorm A \nnorm_{{2}}$ for any symmetric matrix $A\in
{\mathbb{R}}^{a \times a}$.

\section{Background and problem formulation}
\label{SecBackground}

We begin by providing some background on Markov random fields,
defining the problem of graphical model selection and describing our
method based on neighborhood logistic regression.

\subsection{Pairwise Markov random fields}\label{sec21}

Let $X = (X_1, X_2, \ldots, X_p)$ denote a random vector with each
variable $X_s$ taking values in a corresponding set $\mathcal{X}_s$.
Say we are given an undirected graph $G$ with vertex set $V=
\{1,\ldots, p\}$ and edge set $E$, so that each random variable $X_s$
is associated with a vertex $s \in V$. The pairwise Markov random field
associated with the graph $G$ over the random vector $X$ is the family
of distributions of $X$ which factorize as $\mathbb{P}(x) \propto\exp\{
\sum_{(s,t) \in E} \phi_{st}(x_s, x_t) \}$ where for each edge $(s,t)
\in E$, $\phi_{st}$ is a mapping from pairs $(x_s, x_t)
\in\mathcal{X}_s \times \mathcal{X}_t$ to the real line. For models
involving discrete random variables, the pairwise assumption involves
no loss of generality since any Markov random field with higher-order
interactions can be converted (by introducing additional variables) to
an equivalent Markov random field with purely pairwise interactions
(see Wainwright and Jordan \cite{WaiJor03Monster} for details of this
procedure).

\subsubsection*{Ising model}

In this paper, we focus on the special case
of the Ising model in which $X_s \in\{-1,1\}$ for each vertex $s \in
V$, and $\phi_{st}(x_s, x_t) = \theta^*_{st} x_s x_t$ for some
parameter $\theta^*_{st} \in{\mathbb{R}}$, so that the distribution takes
the form
%
%
\begin{equation}
\label{EqnIsing}
\mathbb{P}_{{{\theta^{*}}}}(x) =
\frac{1}{Z({{\theta^{*}}})} \exp \biggl\{\sum _{(s,t)\in E} \theta^*_{st} x_s
x_t \biggr\}.
\end{equation}
The partition function $Z({{\theta^{*}}})$ ensures that the distribution
sums to one. This model is used in many applications
of spatial statistics such as modeling the behavior of gases or
magnets in statistical physics \cite{Ising25}, building statistical
models in computer vision \cite{Geman84} and social network
analysis.

\subsection{Graphical model selection}\label{sec22}

Suppose that we are given a collection $\mathfrak{X}_1^n:=\{x^{(1)},
\ldots, x^{(n)} \}$ of $n$ samples where each
$p$-dimensional vector $x^{(i)} \in\{-1,+1\}^p$ is drawn
in an i.i.d. manner from a distribution $\mathbb{P}_{{{\theta^{*}}}}$
of the
form (\ref{EqnIsing}) for parameter vector ${{\theta^{*}}}$ and graph
$G= (V, E)$ over the $p$ variables. It is
convenient to view the parameter vector ${{\theta^{*}}}$ as a
${p\choose2}$-dimensional vector, indexed by pairs of distinct
vertices but nonzero if and only if the vertex pair $(s,t)$ belongs
to the unknown edge set $E$ of the underlying graph $G$. The
goal of \textit{graphical model selection} is to infer the edge set
$E$. In this paper, we study the
slightly stronger criterion of \textit{signed edge recovery}; in
particular, given a graphical model with parameter ${{\theta^{*}}}$, we
define the edge sign vector
%
%
\begin{equation}
\label{EqnSignedEdge}
E^* :=
\cases{\operatorname{sign}(\theta^*_{st}), &\quad if $(s,t)\in E$, \cr
0, &\quad otherwise.}
\end{equation}
Here the sign function takes value $+1$ if $\theta^*_{st} > 0$, value
$-1$ if $\theta^*_{st} < 0$ and $0$, otherwise. Note that the weaker
graphical model selection problem amounts to recovering the vector
$|E^*|$ of absolute values.

The classical notion of statistical consistency applies to the
limiting behavior of an estimation procedure as the sample size
$n$ goes to infinity with the model size $p$ itself
remaining fixed. In many contemporary applications of graphical
models---among them gene microarray data and social network
analysis---the model dimension $p$ is comparable to or larger than
the sample size $n$, so that the relevance of such ``fixed
$p$'' asymptotics is limited. With this motivation, our analysis
in this paper is of the high-dimensional nature, in which both the
model dimension and the sample size are allowed to increase, and we
study the scalings under which consistent model selection is
achievable.

More precisely, we consider sequences of graphical model selection
problems, indexed by the sample size $n$, number of vertices
$p$ and maximum node degree $d$. We assume that the sample
size $n$ goes to infinity, and both the problem dimension $p
= p(n)$ and $d= d(n)$ may also scale
as a
function of $n$. The setting of fixed $p$ or $d$ is
covered as a special case. Let $\widehat{E}_n$ be an estimator
of the signed edge pattern $E^*$ based on the $n$ samples.
Our goal is to establish sufficient conditions on the scaling of the
triple $(n, p, d)$ such that our proposed estimator is
consistent in the sense that
\[
\mathbb{P} [\widehat{E}_n= E^* ]
\rightarrow
1 \qquad\mbox{as $n\rightarrow+\infty$}.
\]
We sometimes call this property \textit{sparsistency}, as a shorthand for
consistency of the sparsity pattern of the parameter $\theta^*$.

\subsection{Neighborhood-based logistic regression}\label{sec23}

Recovering the signed edge vector $E^*$ of an
undirected graph $G$ is equivalent to recovering,
for each vertex $r\in V$, its \textit{neighborhood set} $\mathcal{N}(r)
:=\{t \in V \mid(r,t) \in E\}$ along with the correct signs
$\operatorname{sign}(\theta^*_{rt})$ for all $t \in\mathcal {N}(r)$. To
capture both the neighborhood structure and sign pattern, we define the
product set of ``signed vertices'' as $\{-1,1\} \times V$. We use the
shorthand ``$\iota r$'' for elements $(\iota, r) \in\{-1,1\} \times V$.
We then define the \textit{signed neighborhood set} as
%
%
\begin{equation}
\label{EqnSignedNeigh}
\mathcal{N}_{\pm}(r) := \{ \operatorname{sign}(\theta
^*_{rt}) t
\mid t \in\mathcal{N}(r) \}.
\end{equation}
Here the sign function has an unambiguous definition, since
$\theta^*_{rt} \neq0$ for all $t \in\mathcal{N}(r)$. Observe
that this signed neighborhood set $\mathcal{N}_{\pm}(r)$ can be recovered
from the sign-sparsity pattern of the $(p-1)$-dimensional
subvector of parameters
\[
\theta^*_{\setminus r} := \{ \theta^*_{ru}, u
\in
V\setminus r \},
\]
associated with vertex $r$. In order to estimate this vector
$\theta^*_{\setminus r}$, we consider the structure of the conditional
distribution of $X_r$ given the other variables
$X_{\setminus r} = \{X_t \mid t \in V\setminus\{r\}
\}$. A simple calculation shows that under the
model (\ref{EqnIsing}), this conditional distribution takes the form
%
%
\begin{equation}
\label{EqnCondDist}
\mathbb{P}_{{{\theta^{*}}}}(x_r \mid x_{\setminus r}) = \frac{
\exp( 2 x_r\sum_{t \in V\setminus r}
\theta^*_{rt} x_t )}{ \exp(2 x_r
\sum_{t \in V\setminus r} \theta^*_{rt} x_t ) + 1}.
\end{equation}
Thus the variable $X_r$ can be viewed as the response variable
in a logistic regression in which all of the other variables
$X_{\setminus r}$ play the role of the covariates.

With this set-up, our method for estimating the sign-sparsity pattern
of the regression vector $\theta^*_{\setminus r}$ and
hence\vspace*{1pt} the neighborhood structure $\mathcal{N}_{\pm}(r)$ is
based on computing an $\ell_1$-regularized logistic regression of $X_r$
on the other variables $X_{\setminus r}$. Explicitly, given
$\mathfrak{X}_1^n= \{ x^{(1)}, x^{(2)}, \ldots, x^{(n)} \}$, a set of
$n$ i.i.d. samples, this regularized regression problem is a convex
program of the form
%
%
\begin{equation}\label{eq:lr}
\min_{\theta_{\setminus r} \in{\mathbb{R}}^{p-1} }
\bigl\{ \ell(\theta; \mathfrak{X}_1^n) + \lambda_{(n,p
,d)}
\|\theta_{\setminus r}\|_1 \bigr\},
\end{equation}
where
%
%
\begin{equation}
\ell(\theta; \mathfrak{X}_1^n) := -\frac{1}{n}
\sum_{i=1}^n\log\mathbb{P}_{\theta}\bigl(x^{(i)}_r \mid
x^{(i)}_{\setminus r}\bigr)
\end{equation}
is the rescaled negative log likelihood (the rescaling factor
$1/n$ in this definition is for later theoretical convenience)
and $\lambda_{(n,p,d)} > 0$ is a regularization
parameter, to be
specified by the user. For notational convenience, we will also use
$\lambda_{n}$
as notation for this regularization parameter suppressing the
potential dependence on $p$ and $d$.

Following some algebraic manipulation, the regularized negative log
likelihood can be written as
%
%
\begin{equation}
\label{EqnRegLikeTwo}
\min_{\theta_{\setminus r} \in{\mathbb{R}}^{p-1} }
\Biggl\{
\frac{1}{n} \sum_{i=1}^nf\bigl(\theta; x^{(i)}\bigr) -
\sum_{u \in V\setminus r} \theta_{ru} \widehat{\mu}_{ru}
+ \lambda_n\|\theta_{\setminus r}\|_1 \Biggr\},
\end{equation}
where
%
%
\begin{equation}
\label{EqnDefnLogFun}
f(\theta; x) := \log\biggl\{\exp\biggl( \sum_{t \in
V\setminus r} \theta_{rt} x_t \biggr) + \exp\biggl(-
\sum_{t \in
V\setminus r} \theta_{rt} x_t \biggr) \biggr\}
\end{equation}
is a rescaled logistic loss, and $\widehat{\mu}_{ru} :=
\frac{1}{n} \sum_{i=1}^nx^{(i)}_rx^{(i)}_u$
are empirical moments. Note the objective
function (\ref{EqnRegLikeTwo}) is convex but not differentiable, due
to the presence of the $\ell_1$-regularizer. By Lagrangian duality,
the problem (\ref{EqnRegLikeTwo}) can be re-cast as a constrained
problem over the ball $\|\theta_{\setminus r}\|_1 \leq
C(\lambda_n)$. Consequently, by the Weierstrass theorem, the
minimum over $\theta_{\setminus s}$ is always achieved.

Accordingly, let $\widehat{\theta}^n_{\setminus r}$ be an element of the
minimizing set of problem (\ref{EqnRegLikeTwo}). Although
$\widehat{\theta}^n_{\setminus r}$ need not be unique in general
since the
problem (\ref{EqnRegLikeTwo}) need not be strictly convex, our
analysis shows that in the regime of interest, this minimizer
$\widehat{\theta}^n_{\setminus r}$ is indeed unique. We use
$\widehat{\theta}^n_{\setminus r}$ to estimate the signed neighborhood
$\mathcal{N}_{\pm}(r)$ according to
%
%
\begin{equation}
\label{EqnNeighEst}
\widehat{\mathcal{N}}_\pm(r) := \{ \operatorname
{sign}(\widehat{\theta}_{ru})
u \mid u \in V\setminus r, \widehat{\theta}_{su} \neq0 \}.
\end{equation}
We say that the full graph $G$ is estimated consistently, written
as the event $\{ \widehat{E}_n= E^* \}$, if every signed
neighborhood is recovered---that is, $\widehat{\mathcal{N}}_\pm(r) =
\mathcal{N}_{\pm}
(r)$
for all $r\in V$.

\section{Method and theoretical guarantees}
\label{SecOutline}

Our main result concerns conditions on the sample size $n$
relative to the parameters of the graphical model---more specifically,
the number of nodes $p$ and maximum node degree $d$---that
ensure that the collection of signed neighborhood
estimates (\ref{EqnNeighEst}), one for each node $r$ of the
graph, agree with the true neighborhoods so that the full graph is
estimated consistently. In this section, we begin by stating the
assumptions that underlie our analysis, and then give a precise
statement of the main result. We then provide a high-level overview of
the key steps involved in its proof, deferring details to later
sections. Our analysis proceeds by first establishing sufficient
conditions for correct signed neighborhood recovery---that is, $\{
\widehat{\mathcal{N}}_\pm(r) = \mathcal{N}_{\pm}(r) \}$---for some
fixed node
$r
\in V$. By showing that this neighborhood consistency is
achieved at sufficiently fast rates, we can then use a union bound
over all $p$ nodes of the graph to conclude that consistent graph
selection is also achieved.

\subsection{Assumptions}\label{sec31}

Success of our method requires certain assumptions on the structure of
the logistic regression problem. These assumptions are stated in
terms of the Hessian of the likelihood function ${\mathbb{E}}\{\log
\mathbb{P}_{\theta}[X_r \mid X_{\setminus r}]\}$ as
evaluated at the true model parameter $\theta^*_{\setminus r}
\in
{\mathbb{R}}^{p-1}$. More specifically, for any fixed node $r\in
V$, this Hessian is a $(p-1) \times(p-1)$ matrix of the
form
%
%
\begin{equation}
\label{EqnDefnQstar}
Q^{*}_r := {\mathbb{E}}_{{{\theta^{*}}}} \{ \nabla^2 \log
\mathbb{P}_{{{\theta^{*}}}}[X_r \mid X_{\setminus r}] \}.
\end{equation}
For future reference, this is given as the explicit expression
%
%
\begin{equation}
\label{EqnDefnQstarExp}
Q^{*}_r = {\mathbb{E}}_{{{\theta^{*}}}} [\eta(X; {{\theta^{*}}})
X_{\setminus r} X_{\setminus r}^T ],
\end{equation}
where
%
%
\begin{equation}
\label{EqnDefnVarfun}
\eta(u; \theta) := \frac{ 4 \exp(2 u_r
\sum_{t \in V\setminus r} \theta_{rt} u_t
)}{ (\exp(2 u_r\sum_{t \in V\setminus r}
\theta_{rt} u_t ) + 1 )^2 }
\end{equation}
is the variance function. Note that the matrix $Q^{*}_r$ is the
Fisher information matrix associated with the local conditional
probability distribution. Intuitively, it serves as the counterpart
for discrete graphical models of the covariance matrix ${\mathbb{E}}[X X^T]$
of Gaussian graphical models, and indeed our assumptions are analogous
to those imposed in previous work on the Lasso for Gaussian linear
regression \cite{Meinshausen06,Tropp06,ZhaoYu06}.

In the following we write simply $Q^{*}$ for the matrix $Q^{*}
_r$ where
the reference node $r$ should be understood implicitly.
Moreover, we use $S:=\{ (r, t) \mid t \in
\mathcal{N}(r) \}$ to denote the subset of indices associated with edges
of $r$, and ${{S^{c}}}$ to denote its complement. We use
$Q^{*}_{SS}$ to denote the $d\times d$
sub-matrix of $Q^{*}$ indexed by $S$. With this notation, we
state our assumptions:

\subsubsection*{(\textup{A}1) Dependency condition}

The subset of the Fisher information matrix corresponding to the
relevant covariates has bounded eigenvalues; that is, there exists a
constant $C_{\min}> 0$ such that
%
%
\begin{equation}
\label{EqnEigBounds}
\Lambda_{\min}(Q^{*}_{SS}) \geq C_{\min}.
\end{equation}
Moreover, we require that
$\Lambda_{\max}({\mathbb{E}}_{{{\theta^{*}}}}[X_{\setminus r}
X_{\setminus r}^T])
\leq D_{\max}$. These conditions ensure that the relevant covariates do
not become overly dependent. (As stated earlier, we have suppressed
notational dependence on $r$; thus these conditions are assumed
to hold for each $r\in V$.)

\subsubsection*{(\textup{A}2) Incoherence condition}

Our next assumption captures the intuition that the large number of
irrelevant covariates (i.e., nonneighbors of node $r$) cannot exert an
overly strong effect on the subset of relevant covariates (i.e.,
neighbors of node $r$). To formalize this intuition, we require the
existence of an $\alpha\in(0,1]$ such that
%
%
\begin{equation}
\label{EqnDefnMutualInco}
\nnorm Q^*_{{{S^{c}}}S} (Q^*_{S S})^{-1} \nnorm_{{\infty}} \leq1
- \alpha.
\end{equation}
%

\subsection{Statement of main result}\label{sec32}

We are now ready to state our main result on the performance of
neighborhood logistic regression for graphical model selection.
Naturally, the limits of model selection are determined by the minimum
value over the parameters $\theta^*_{rt}$ for pairs $(r,
t)$ included in the edge set of the true graph. Accordingly, we
define the parameter
%
%
\begin{equation}
\label{EqnDefnTmin}
\theta^*_{\min} = {\min_{(r, t) \in E}} |\theta^*_{rt}|.
\end{equation}
With this definition, we have the following:
\begin{theos}
\label{ThmMain}
Consider an Ising graphical model with parameter vector $\theta^*$
and associated edge set $E^*$ such that conditions \textup{(A1)} and
\textup{(A2)} are satisfied by the population Fisher information matrix
$Q^{*}$, and let $\mathfrak{X}_1^n$ be a set of $n$ i.i.d. samples from
the model specified by $\theta^*$. Suppose that the regularization
parameter $\lambda_n$ is selected to satisfy
%
%
\begin{equation}\label{eq:reg}
\lambda_n \geq\frac{16
(2-\alpha)}{\alpha} \sqrt{\frac{\log p}{n}}.
\end{equation}
Then there exist positive constants
$L$ and $K$, independent of $(n, p, d)$, such
that if
%
%
\begin{equation}
\label{EqnGrowthCondition}
n> L d^3 \log p,
\end{equation}
then the following properties hold with probability at least $1 - 2
\exp(-K\lambda_{n}^2 n )$.
\begin{enumerate}[(a)]
\item[(a)] For each node $r\in V$, the $\ell_1$-regularized
logistic regression (\ref{eq:lr}), given data~$\mathfrak{X}_1^n$, has
a unique solution, and so uniquely
specifies a signed neighborhood $\widehat{\mathcal{N}}_\pm(r)$.
\item[(b)] For each $r\in V$, the estimated signed
neighborhood $\widehat{\mathcal{N}}_\pm(r)$ correctly excludes all edges
\textit{not} in the true neighborhood. Moreover, it correctly
includes all edges $(r,t)$ for which
$|\theta_{rt}^*| \geq\frac{10}{C_{\min}} \sqrt{d}
\lambda_{n}$.
\end{enumerate}
\end{theos}

The theorem not only specifies sufficient conditions but also the
probability with which the method recovers the true signed
edge-set. This probability decays exponentially as a function of
$\lambda_{n}^2 n$ which leads naturally to the following
corollary on model selection consistency of the method for a sequence
of Ising models specified by $(n, p(n),
d(n))$.
\begin{cors}
Consider a sequence of Ising models with graph edge sets
$\{E^{*}_{p(n)} \}$ and parameters $ \{ \theta^*_{(n,p, d)} \}$; each of
which satisfies conditions \textup{(A1)} and~\textup{(A2)}. For each
$n$, let $\mathfrak{X}_1^n$ be a set of $n$ i.i.d. samples from the
model specified by $\theta^*_{(n,p, d)}$, and suppose that $(n,
p(n),d(n))$ satisfies the scaling condition (\ref{EqnGrowthCondition})
of Theorem \ref{ThmMain}. Suppose further that the sequence
$\{\lambda_n\}$ of regularization parameters satisfies condition
(\ref{eq:reg}) and
%
%
\begin{equation}\label{eq:regcond}
\lambda_{n}^2 n \rightarrow\infty
\end{equation}
and the minimum parameter weights satisfy
%
%
\begin{equation}\label{eq:thetacond}
{\min_{(r,t)\in E^{*}_{n}}} \bigl|\theta^*_{(n,p, d)}(r,t)\bigr| \geq
\frac{10}{C_{\min}} \sqrt{d} \lambda_{n}
\end{equation}
for sufficiently large $n$. Then the method is model selection
consistent so that if $\widehat{E}_{p(n)}$ is the graph
structure estimated
by the method given data $\mathfrak{X}_1^n$, then $\mathbb{P}
[\widehat{E }_{p(n)} =
E^*_{p(n)} ] \rightarrow1$ as $n
\rightarrow\infty$.
\end{cors}

\subsubsection*{Remarks}

(a) It is worth noting that the scaling
condition (\ref{EqnGrowthCondition}) on $(n, p, d)$
allows for graphs and sample sizes in the ``large $p$, small
$n$'' regime (meaning $p\gg n$), as long as the degrees
are bounded, or grow at a sufficiently slow rate. In particular, one
set of sufficient conditions are the scalings
\[
d= O(n^{c_1}) \quad\mbox{and}\quad p=
O(e^{n^{c_2}}),\qquad 3 c_1+ c_2< 1,
\]
for some constants $c_1,c_2> 0$. Under these scalings, note
that we have $d^3 \log(p) = \mathcal{O}(n^{3
c_1+ c_2}) = o(n)$, so that
condition (\ref{EqnGrowthCondition}) holds.

A bit more generally, note that in the regime $p\gg n$, the
growth condition (\ref{EqnGrowthCondition}) requires that that
$d= o(p)$. However, in many practical applications of
graphical models (e.g., image analysis, social networks), one is
interested in node degrees $d$ that remain bounded or grow
sub-linearly in the graph size so that this condition is not
unreasonable.

(b) Loosely stated, the theorem requires that the edge weights
are not too close to zero (in absolute value) for the method to
estimate the true graph. In particular, conditions (\ref{eq:reg}) and
(\ref{eq:thetacond}) imply that the minimum edge weight $\theta
^*_{\min}$ is
required to scale as
\[
\theta^*_{\min}= \Omega\Biggl(\sqrt{\frac{d \log
p}{n}} \Biggr).
\]
Note that in the classical fixed $(p,d)$ case, this
reduces to the familiar scaling requirement of $\theta^*_{\min}=
\Omega(n^{-1/2})$.

(c) In the high-dimensional setting (for $p\rightarrow
+\infty)$, a choice of the regularization parameter satisfying both
conditions (\ref{eq:reg}) and (\ref{eq:regcond}) is, for example,
\[
\lambda_{n} = \frac{16 (2-\alpha)}{\alpha} \sqrt{\frac
{\log p}{n}}
\]
for which the probability of incorrect model selection decays at
rate\break
${\mathcal{O}}(\exp(- K' \log p))$ for some constant $K' > 0$.
In the classical setting (fixed~$p$), this choice can be modified
to $\lambda_{n} = \frac{16 (2-\alpha)}{\alpha}
\sqrt{\frac{\log( pn)}{n}}$.

The analysis required to prove Theorem \ref{ThmMain} can be divided
naturally into two parts. First, in Section \ref{SecFixedDesign}, we
prove a result (stated as Proposition \ref{PropFixed}) for ``fixed
design'' matrices. More precisely, we show that if the dependence
condition (A1) and the mutual incoherence condition \textup{(A2)}
hold for the \textit{sample Fisher information matrix}
%
%
\begin{equation}
\label{EqnDefnSampleFisher}
Q^{n} := \widehat{{\mathbb{E}}} [\eta(X; {{\theta^{*}}})
X_{\setminus r} X_{\setminus r}^T ] =
\frac{1}{n} \sum_{i=1}^n\eta\bigl(x^{(i)}; {{\theta^{*}}}\bigr)
x^{(i)}_{\setminus r} \bigl(x^{(i)}_{\setminus r}\bigr)^T,
\end{equation}
then the growth condition (\ref{EqnGrowthCondition}) and choice of
$\lambda_n$ from Theorem \ref{ThmMain} are sufficient to ensure
that the graph is recovered with high probability.

The second part of the analysis, provided in
Section \ref{SecPopulation}, is devoted to showing that under the
specified growth condition (\ref{EqnGrowthCondition}), imposing
incoherence and
dependence assumptions on the \textit{population version} of the Fisher
information $Q^{*}$ guarantees (with high probability) that analogous
conditions hold for the sample quantities~$Q^{n}$. On one hand, it
follows immediately from the law of large numbers that the empirical
Fisher information $Q^{n}_{AA}$ converges to the population version
$Q^{*}_{AA}$ for any \textit{fixed} subset $A$. However, in the
current setting, the added delicacy is that we are required to control
this convergence over subsets of increasing size. Our proof therefore
requires some large-deviation analysis for random matrices with
dependent elements so as to provide exponential control on the rates
of convergence.

\subsection{Primal-dual witness for graph recovery}
\label{SecPrimalDual}

At the core of our proof lies the notion of a primal-dual witness
used in previous work on the Lasso \cite{Wainwright06new}. In
particular, our proof involves the explicit construction of an optimal
\textit{primal-dual pair}---namely, a primal solution $\widehat{\theta}
\in{\mathbb{R}}^{p- 1}$ along with an associated subgradient vector
$\widehat{z}\in{\mathbb{R}}^{p-1 }$ (which can be interpreted as a dual
solution), such that the sub-gradient optimality conditions associated
with the convex program (\ref{EqnRegLikeTwo}) are satisfied.
Moreover, we show that under the stated assumptions on $(n,
p, d)$, the primal-dual pair $(\widehat{\theta}, \widehat{z})$ can
be constructed such that they act as a \textit{witness}---that is, a
certificate guaranteeing that the method correctly recovers the graph
structure.

For the convex program (\ref{EqnRegLikeTwo}), the zero sub-gradient
optimality conditions \cite{Rockafellar} take the form
%
%
\begin{equation}
\label{EqnZeroGrad}
\nabla\ell(\widehat{\theta}) + \lambda_n\widehat{z} = 0,
\end{equation}
where the dual or subgradient vector $\widehat{z}\in{\mathbb{R}}^{p-1}$
must satisfy the properties
%
%
\begin{equation}
\label{EqnDualvecProp}
\widehat{z}_{rt} = \operatorname{sign}(\widehat{\theta}_{rt})
\qquad\mbox{if
$\widehat{\theta}_i \neq0$}\quad \mbox{and}\quad |\widehat{z}_{rt}|
\leq1\qquad
\mbox{otherwise}.
\end{equation}
By convexity, a pair $(\widehat{\theta}, \widehat{z}) \in{\mathbb{R}}^{p-1}
\times{\mathbb{R}}^{p-1}$ is a primal-dual optimal solution to the
convex program and its dual if and only if the two
conditions (\ref{EqnZeroGrad}) and (\ref{EqnDualvecProp}) are
satisfied. Of primary interest to us is the property that such an
optimal primal-dual pair correctly specifies the signed neighborhood of node
$r$; the necessary and sufficient conditions for such correctness
are
%
%
\begin{subequation}
\begin{eqnarray}
\label{EqnCorrectSign}
\operatorname{sign}(\widehat{z}_{rt}) & = & \operatorname
{sign}(\theta^*_{rt})\qquad
\forall(r, t) \in S:=\{(r, t) \in E\}\quad
\mbox{and} \\
\label{EqnNoFalsePos}
\widehat{\theta}_{ru} & = & 0 \qquad\mbox{for all $(r, u) \in
{{S^{c}}}:=E\setminus S$.}
\end{eqnarray}
\end{subequation}

The $\ell_1$-regularized logistic regression
problem (\ref{EqnRegLikeTwo}) is convex; however, for $p\gg
n$, it need not be strictly convex, so that there may be
multiple optimal solutions. The following lemma, proved in
Appendix \ref{AppLemUnique}, provides sufficient conditions for shared sparsity
among optimal solutions, as well as uniqueness of the optimal
solution:
\begin{lems}
\label{LemUnique}
Suppose that there exists an optimal primal solution $\widehat{\theta}$
with associated optimal dual vector $\widehat{z}$ such that
$\|\widehat{z}_{{{S^{c}}}}\|_\infty< 1$. Then any optimal primal
solution $\widetilde{\theta}$ must have $\widetilde{\theta
}_{{{S^{c}}}} = 0$.
Moreover, if the Hessian sub-matrix $[\nabla^2
\ell(\widehat{\theta})]_{SS}$ is strictly positive
definite, then
$\widehat{\theta}$ is the unique optimal solution.
\end{lems}

Based on this lemma, we construct a primal-dual witness
$(\widehat{\theta},\widehat{z})$ with the following steps.
\begin{enumerate}[(a)]
\item[(a)] First, we set $\widehat{\theta}_{S}$ as the minimizer of the
partial penalized likelihood
%
%
\begin{equation}
\label{thetaSCon}
\widehat{\theta}_{S} = \mathop{\arg\min}_{(\theta_{S}, 0) \in{\mathbb{R}}
^{p-1}} \{\ell(\theta; \mathfrak{X}_1^n) +
\lambda_n\|\theta_{S}\|_1 \}
\end{equation}
and set $\widehat{z}_{S} = \operatorname{sign}(\widehat{\theta}_{S})$.
\item[(b)] Second, we set $\widehat{\theta}_{{S^{c}}}= 0$ so that
condition (\ref{EqnNoFalsePos}) holds.
\item[(c)] In the third step, we obtain $\widehat{z}_{{{S^{c}}}}$ from
(\ref{EqnZeroGrad}) by substituting in the values of
$\widehat{\theta}$ and $\widehat{z}_{S}$. Thus our
construction satisfies conditions (\ref{EqnNoFalsePos}) and
(\ref{EqnZeroGrad}).
\item[(d)] The final and most challenging step consists of showing
that the stated scalings of $(n, p, d)$ imply that,
with high-probability, the remaining conditions (\ref{EqnCorrectSign})
and (\ref{EqnDualvecProp}) are satisfied.
\end{enumerate}
Our analysis in step (d) guarantees that
$\|\widehat{z}_{S^c}\|_\infty< 1$ with high probability. Moreover,
under the conditions of Theorem \ref{ThmMain}, we prove that the
sub-matrix of the sample Fisher information matrix is strictly
positive definite with high probability so that by
Lemma \ref{LemUnique}, the primal solution $\widehat{\theta}$ is
guaranteed to be unique.

It should be noted that, since $S$ is unknown, the primal-dual witness
method is \textit{not} a practical algorithm that could ever be
implemented to solve $\ell_1$-regularized logistic regression.
Rather, it is a proof technique that allows us to establish sign
correctness of the unique optimal solution.

\section{Analysis under sample Fisher matrix assumptions}
\label{SecFixedDesign}

We begin by establishing model selection consistency when assumptions
are imposed directly on the sample Fisher matrix $Q^{n}$, as opposed
to on the population matrix $Q^{*}$, as in Theorem \ref{ThmMain}. In
particular, recalling the definition (\ref{EqnDefnSampleFisher}) of
the sample Fisher information matrix $Q^{n}= \widehat{{\mathbb
{E}}}[\nabla^2
\ell({{\theta^{*}}})]$, we define the ``good event,''
%
%
\begin{equation}
\mathcal{M}(\mathfrak{X}_1^n) := \bigl\{ \mathfrak{X}_1^n\in\{
-1,+1\}^{n\times
p} \mid\mbox{$Q^{n}$ satisfies (A1) and (A2)} \bigr\}.
\end{equation}
As in the statement of Theorem \ref{ThmMain}, the quantities $L$ and
$K$ refer to constants independent of $(n, p,
d)$. With this notation, we have the following:
\begin{props}[(Fixed design)]
\label{PropFixed}
If the event $\mathcal{M}(\mathfrak{X}_1^n)$ holds, the sample size satisfies
$n> L d^2 \log(p)$, and the regularization parameter
is chosen such that $\lambda_n\geq\frac{16
(2-\alpha)}{\alpha} \sqrt{\frac{\log p}{n}}$, then with
probability at least $1 - 2 \exp(- K\lambda_n^2 n
) \rightarrow1$, the following properties hold.

\textup{(a)} For\vspace*{1pt} each node $r\in V$, the $\ell_1$-regularized
logistic regression has a unique solution, and so uniquely
specifies a signed neighborhood $\widehat{\mathcal{N}}_\pm(r)$.

\textup{(b)} For each $r\in V$, the estimated signed
neighborhood vector $\widehat{\mathcal{N}}_\pm(r)$ correctly
excludes all
edges \textit{not} in the true neighborhood. Moreover, it correctly
includes all edges with $|\theta_{rt}| \geq
\frac{10}{C_{\min}} \sqrt{d} \lambda_n$.
\end{props}

Loosely stated, this result guarantees that if the sample Fisher
information matrix is ``good,'' then the conditional probability of
successful graph recovery converges to zero at the specified rate.
The remainder of this section is devoted to the proof of
Proposition \ref{PropFixed}.

\subsection{Key technical results}\label{sec41}
We begin with statements of some key technical lemmas that are central to
our main argument with their proofs deferred to
Appendix \ref{AppTechnical}. The central object is the following
expansion obtained by re-writing the zero-subgradient condition as
%
%
\begin{equation}
\label{EqnTemp}
\nabla\ell(\widehat{\theta}; \mathfrak{X}_1^n) - \nabla\ell
({{\theta^{*}}}; \mathfrak{X}_1^n
) = W^n- \lambda_n\widehat{z},
\end{equation}
where we have introduced the short-hand notation $W^n= -\nabla
\ell({{\theta^{*}}}; \mathfrak{X}_1^n)$ for the
$(p-1)$-dimensional score
function,
\[
W^n := -\frac{1}{n} \sum_{i=1}^n
x^{(i)}_{\setminus r} \biggl\{x^{(i)}_r- \frac{\exp
(\sum_{t
\in V\setminus r} \theta^*_{rt} x^{(i)}_t )
-\exp(-\sum_{t \in V\setminus r} \theta^*_{rt}
x^{(i)}_t ) } {\exp(\sum_{t \in V\setminus r}
\theta^*_{rt} x^{(i)}_t ) +\exp(-\sum_{t \in
V\setminus
r} \theta^*_{rt} x^{(i)}_t ) } \biggr\}.
%
\]
For future reference, note that ${\mathbb{E}}_{{{\theta^{*}}}}[W^n] = 0$.
Next, applying the mean-value theorem coordinate-wise to the
expansion (\ref{EqnTemp}) yields
%
%
\begin{equation}
\label{EqnTay}
\nabla^2 \ell({{\theta^{*}}}; \mathfrak{X}_1^n) [\widehat
{\theta} - {{\theta^{*}}}
] = W^n-\lambda_n\widehat{z}+ R^n,
\end{equation}
where the remainder term takes the form
%
%
\begin{equation}
\label{EqnRemainder}
R^n_j = \bigl[\nabla^2 \ell\bigl(\onebar{\theta}^{(j)};
\mathfrak{X}_1^n\bigr) -
\nabla^2 \ell({{\theta^{*}}}; \mathfrak{X}_1^n) \bigr]_j^T
(\widehat{\theta}-
{{\theta^{*}}})
\end{equation}
with $\onebar{\theta}^{(j)}$ a parameter vector on the line between
${{\theta^{*}}}$
and $\widehat{\theta}$, and with $[\cdot]_j^T$ denoting the $j$th
row of
the matrix.
The following lemma addresses the behavior of the
term $W^n$ in this expansion:
\begin{lems}
\label{LemBernOne}
For the specified mutual incoherence parameter $\alpha\in(0,1]$, we
have
%
%
\begin{equation}
\mathbb{P} \biggl(\frac{2-\alpha}{\lambda_n} \|W^n\|
_\infty\geq
\frac{\alpha}{4} \biggr) \leq2 \exp\biggl( - \frac{\alpha^2
\lambda_n^2}{128 (2-\alpha)^2} n+ \log(p) \biggr),
\end{equation}
which converges to zero at rate $\exp(-c\lambda_n^2
n)$ as long as $\lambda_n\geq\frac{16
(2-\alpha)}{\alpha} \sqrt{\frac{\log p}{n}}$.
\end{lems}

See Appendix \ref{AppBernOne} for the proof of this claim.

The following lemma establishes that the sub-vector
$\widehat{\theta}_S$ is an $\ell_2$-consistent estimate of the
true sub-vector $\theta^*_S$:
\begin{lems}[($\ell_2$-consistency of primal subvector)]
\label{Leml2cons}
If $\lambda_nd\leq\frac{C_{\min}^2}{10 D_{\max}}$ and\break
$\|W^n\|_\infty\leq\lambda_n/4$, then
%
%
\begin{equation}
\label{Eqnl2cons}
\|\widehat{\theta}_{S} - \theta_{S}\|_{2} \leq
\frac{5}{C_{\min}} \sqrt{d} \lambda_{n}.
\end{equation}
\end{lems}

See Appendix \ref{Appl2cons} for the proof of this claim.

Our final technical lemma provides control on the remainder
term (\ref{EqnRemainder}).
\begin{lems}
\label{LemTayRem}
If $\lambda_nd\leq\frac{C_{\min}^2}{100 D_{\max}}
\frac{\alpha}{2-\alpha}$ and $\|W^n\|_\infty\leq
\lambda_n/4$, then
\[
\frac{\|R^n\|_\infty}{\lambda_n} \leq\frac{25
D_{\max}}{C_{\min}^2} \lambda_nd \leq\frac{\alpha}{4
(2-\alpha)}.
\]
\end{lems}

See Appendix \ref{AppTayRem} for the proof of this claim.


\subsection[Proof of Proposition 1]{Proof of Proposition \protect\ref{PropFixed}}\label{sec42}

Using these lemmas, the proof of Proposition~\ref{PropFixed} is
straightforward. Consider the choice of the regularization parameter,
$\lambda_n= 16 \frac{2-\alpha}{\alpha} \sqrt{\frac{\log
p}{n}}$. This choice satisfies the condition of
Lemma \ref{LemBernOne}, so that we may conclude that with probability
greater than $1-2\exp(- c\lambda_n^2 n)
\rightarrow
1$, we have
\[
\|W^n\|_\infty\leq\frac{\alpha}{2-\alpha}
\frac{\lambda}{4} \leq\frac{\lambda}{4}
\]
using the fact that $\alpha\leq1$. The remaining two conditions
that we need to apply the technical lemmas concern upper bounds on
the quantity $\lambda_nd$. In particular, for a sample size
satisfying $n> \frac{100^2 D_{\max}^2}{C_{\min}^4}
\frac{(2-\alpha)^4}{\alpha^4} d^2 \log p$, we have
\begin{eqnarray*}
\lambda_nd& = & \frac{16(2-\alpha)}{\alpha} \sqrt
{\frac{\log
p}{n}} d \\
%
%
& \leq& \frac{16 C_{\min}^2}{100 D_{\max}} \frac{\alpha
}{(2-\alpha)}
\\
& < & \frac{C_{\min}^2}{10 D_{\max}}
\end{eqnarray*}
so that the conditions of both Lemmas \ref{Leml2cons} and
\ref{LemTayRem} are satisfied.

We can now proceed to the proof of Proposition \ref{PropFixed}.
Recalling our shorthand $Q^{n}= \nabla^2_\theta
\ell({{\theta^{*}}}; \mathfrak{X}_1^n)$ and the fact that we have set
$\widehat{\theta}_{S^{c}}= 0$ in our primal-dual construction, we can
re-write condition (\ref{EqnTay}) in block form as
%
%
\begin{subequation}
\label{EqnBlock}
\begin{eqnarray}
\label{EqnBlockA}
Q^{n}_{{S^{c}}S} [\widehat{\theta}_S-\theta^*_S] & = &
W^n_{{S^{c}}} -\lambda_n\widehat{z}_{S^{c}}+
R^n_{{S^{c}}}, \\
\label{EqnBlockB}
Q^{n}_{SS} [\widehat{\theta}_S-\theta^*_S
] & = & W^n_S-\lambda_n
\widehat{z}_S+ R^n_S.
\end{eqnarray}
\end{subequation}
Since the matrix $Q^{n}_{SS}$ is invertible by assumption,
the conditions (\ref{EqnBlock}) can be re-written as
%
%
\begin{equation}
Q^{n}_{{S^{c}}S} (Q^{n}_{SS})^{-1}
[W^n_S-\lambda_n\widehat{z}_S+
R^n_S ] = W^n_{{S^{c}}} -\lambda
_n
\widehat{z}_{S^{c}}+ R^n_{{S^{c}}}.
\end{equation}
Rearranging yields the condition,
%
%
\begin{equation}
\label{EqnDefnDualsub}\qquad
[W^n_{{S^{c}}} - R^n_{{S^{c}}} ] -
Q^{n}_{{S^{c}}S} (Q^{n}_{SS})^{-1}
[W^n_S- R^n_S ] +
\lambda_nQ^{n}_{{S^{c}}S} (Q^{n}_{SS})^{-1}
\widehat{z}_S= \lambda_n\widehat{z}_{S^{c}}.
\end{equation}

\subsubsection*{Strict dual feasibility}

We now demonstrate that for the
dual sub-vector $\widehat{z}_{S^{c}}$ defined by
(\ref{EqnDefnDualsub}), we have $\|\widehat{z}_{S^{c}}\|
_\infty<
1$. Using the triangle inequality and the mutual incoherence
bound (\ref{EqnDefnMutualInco}), we have that
%
%
\begin{eqnarray}
\label{EqnTriangle}
\|\widehat{z}_{S^{c}}\|_\infty& \leq& \nnorm Q^{n}_{{S^{c}}S}
(Q^{n}_{SS})^{-1} \nnorm_{{\infty}}
\biggl[\frac{\|W^n_S\|_\infty}{\lambda_n} +
\frac{\|R^n_S\|_\infty}{\lambda_n} + 1
\biggr]\\
&&{} + \frac{\|R^n_{S^{c}}\|_\infty}{\lambda_n} +
\frac{\|W^n_{S^{c}}\|_\infty}{\lambda_n} \nonumber\\
& \leq& (1-\alpha) + (2-\alpha)
\biggl[\frac{\|R^n\|_\infty}{\lambda_n} +
\frac{\|W^n\|_\infty}{\lambda_n} \biggr].
\end{eqnarray}
Next, applying Lemmas \ref{LemBernOne} and \ref{LemTayRem}, we have
\[
\|\widehat{z}_{S^{c}}\|_\infty\leq(1-\alpha) + \frac{\alpha
}{4} +
\frac{\alpha}{4} = 1 - \frac{\alpha}{2}
\]
with probability converging to one.

\subsubsection*{Correct sign recovery}

We next show that our primal sub-vector $\widehat{\theta}_S$ defined by
(\ref{thetaSCon}) satisfies sign consistency, meaning that
$\operatorname{sgn}(\widehat{\theta}_{S}) = \operatorname
{sgn}(\theta^*_{S})$. In order to do so, it suffices to show that
\[
\|\theta_{S} - \theta^*_{S} \|_\infty\leq
\frac{\theta^*_{\min}}{2}
\]
recalling the notation $\theta^*_{\min}:=\min_{(r, t) \in
E} |\theta^*_{rt}|$. From Lemma \ref{Leml2cons}, we have
$\|\theta_{S} - \theta^*_{S}\|_{2} \leq\frac{5}{C_{\min}}
\sqrt{d} \lambda_{n}$ so that
\begin{eqnarray*}
\frac{2}{\theta^*_{\min}}\|\theta_{S} - \theta^*_{S}\|_\infty&
\leq
& \frac{2}{\theta^*_{\min}}\|\theta_{S} - \theta^*_{S}\|_2 \\
& \leq& \frac{2}{\theta^*_{\min}} \frac{5}{C_{\min}} \sqrt{d}
\lambda_n,
\end{eqnarray*}
which is less than one as long as $\theta^*_{\min}\geq\frac
{10}{C_{\min}}
\sqrt{d} \lambda_n$.

\section{Uniform convergence of sample information matrices}
\label{SecPopulation}

In this section we complete the proof of Theorem \ref{ThmMain} by
showing that if the dependency (A1) and incoherence
(A2) assumptions are imposed on the \textit{population} Fisher
information matrix then under the specified scaling of $(n,
p, d)$, analogous bounds hold for the \textit{sample} Fisher
information matrices with probability converging to one. These
results are not immediate consequences of classical random matrix
theory (e.g., \cite{DavSza01}) since the elements of $Q^{n}$ are
highly dependent.
Recall the definitions
%
%
\begin{equation}
\label{EqnRecall}
Q^{*}:={\mathbb{E}}_{{{\theta^{*}}}} [\eta(X; {{\theta^{*}}})
X_{\setminus r} X_{\setminus r}^T ]\quad \mbox{and}\quad
Q^{n}:=\widehat{{\mathbb{E}}} [\eta(X; {{\theta^{*}}})
X_{\setminus r} X_{\setminus r}^T ],
\end{equation}
where ${\mathbb{E}}_{{\theta^{*}}}$ denotes the population
expectation, and
$\widehat{{\mathbb{E}}}$ denotes the empirical expectation, and the variance
function $\eta$ was defined previously in (\ref{EqnDefnVarfun}). The
following lemma asserts that the eigenvalue bounds in assumption
(A1) hold with high probability for sample covariance matrices:
\begin{lems}
\label{LemConcentrateEig}
Suppose that assumption \textup{(A1)} holds for the population matrix
$Q^{*}$ and ${\mathbb{E}}_{{{\theta^{*}}}}[X X^T]$. For any $\delta>
0$ and some
fixed constants $A$ and $B$, we have
%
%
\begin{subequation}
\begin{eqnarray}\quad 
\label{EqnEigConcMax}
\mathbb{P} \Biggl[\Lambda_{\max} \Biggl[\frac{1}{n} \sum_{i=1}^n
x^{(i)}_{\setminus r} \bigl(x^{(i)}_{\setminus r}\bigr)^T \Biggr]
\geq D_{\max}
+ \delta\Biggr] &\leq& 2 \exp\biggl(-A\frac{\delta^2 n}{
d^2}
+ B \log(d) \biggr), \hspace*{-30pt}\\
\label{EqnEigConcMin}
\mathbb{P}[\Lambda_{\min}(Q^{n}_{SS}) \leq C_{\min}- \delta] &
\leq& 2
\exp\biggl(-A\frac{\delta^2 n}{ d^2} + B \log(d)
\biggr).\hspace*{-30pt}
\end{eqnarray}
\end{subequation}
\end{lems}

The following result is the analog for the incoherence assumption
(A2) showing that the scaling of $(n, p, d)$
given in Theorem \ref{ThmMain} guarantees that population incoherence
implies sample incoherence.
\begin{lems}
\label{LemConcentrateInco}
If the population covariance satisfies a mutual incoherence
condition (\ref{EqnDefnMutualInco}) with parameter $\alpha\in(0,1]$
as in assumption \textup{(A2)}, then the sample matrix satisfies an
analogous version, with high probability in the sense that
%
%
\begin{equation}
\label{EqnSampleMutinc}
\mathbb{P} \biggl[\nnorm Q^{n}_{{S^{c}}S} (Q^{n}_{S S})^{-1} \nnorm
_{{\infty}} \geq1- \frac{\alpha}{2} \biggr]
\leq\exp\biggl(-K\frac{n}{d^3} + \log(p) \biggr).
\end{equation}
\end{lems}

Proofs of these two lemmas are provided in the following sections.
Before proceeding, we take note of a simple bound to be
used repeatedly throughout our arguments. By definition of the
matrices $Q^{n}(\theta)$ and $Q(\theta)$ [see
(\ref{EqnDefnSampleFisher}) and (\ref{EqnDefnQstarExp})],
the $(j,k)$th element of the difference matrix $Q^{n}(\theta) -
Q(\theta)$ can be written as an i.i.d. sum of the form $Z_{jk} =
\frac{1}{n}\sum_{i=1}^nZ^{(i)}_{jk}$ where each $Z^{(i)}_{jk}$
is zero-mean and bounded (in particular, $|Z^{(i)}_{jk}| \leq4$). By the
Azuma--Hoeffding bound \cite{Hoeffding63}, for any indices $j,k = 1,
\ldots, d$ and for any $\varepsilon> 0$, we have
%
%
\begin{equation}
\label{EqnAzuma}
\mathbb{P}[(Z_{jk})^2 \geq\varepsilon^2 ] =
\mathbb{P}\Biggl[\Biggl|\frac{1}{n} \sum_{i=1}^nZ^{(i)}_{jk}\Biggr| \geq
\varepsilon\Biggr]
\leq2 \exp\biggl( -\frac{\varepsilon^2 n}{32} \biggr).
\end{equation}
So as to simplify notation, throughout this section, we use $K$
to denote a universal positive constant, independent of $(n,
p, d)$. Note that the precise value and meaning of
$K$ may differ from line to line.

\subsection[Proof of Lemma 5]{Proof of Lemma \protect\ref{LemConcentrateEig}}\label{sec51}

By the Courant--Fischer variational representation~\cite{Horn85}, we
have
\begin{eqnarray*}
\Lambda_{\min}(Q_{SS}) & = & \min_{\|x\|_2 = 1} x^T
Q_{SS} x \\
& = & \min_{\|x\|_2 = 1} \{ x^T Q^{n}_{SS} x +
x^T(Q_{SS} - Q^{n}_{SS})x \} \\
& \leq& y^T Q^{n}_{SS} y + y^T(Q_{SS} -
Q^{n}_{SS}) y,
\end{eqnarray*}
where $y \in{\mathbb{R}}^{d}$ is a unit-norm minimal eigenvector of
$Q^{n}_{SS}$. Therefore, we have
\[
\Lambda_{\min}(Q^{n}_{SS}) \geq\Lambda_{\min}(Q_{S
S})
- \nnorm Q_{SS}- Q^{n}_{SS} \nnorm_{{2}} \geq
C_{\min}- \nnorm Q_{SS}- Q^{n}_{SS} \nnorm_{{2}}.
\]
Hence it suffices to obtain a bound on the spectral norm
$\nnorm Q_{SS}- Q^{n}_{SS} \nnorm_{{2}}$. Observe
that
\[
\nnorm Q^{n}_{SS} - Q_{SS} \nnorm_{{2}} \leq
\Biggl( \sum_{j=1}^d\sum_{k=1}^d(Z_{jk})^2
\Biggr)^{{1/2}}.
\]
Setting $\varepsilon^2 = \delta^2/d^2$ in (\ref{EqnAzuma})
and applying the union bound over the $d^2$ index pairs $(j,k)$ then
yields
%
%
\begin{equation}
\label{EqnUseful}
\mathbb{P}[\nnorm Q^{n}_{SS} -Q_{SS} \nnorm_{{2}}
\geq
\delta] \leq2 \exp\biggl(-K\frac{\delta^2
n}{d^2} + 2 \log(d) \biggr).
\end{equation}

Similarly, we have
\begin{eqnarray*}
&& \mathbb{P}\Biggl[\Lambda_{\max}\Biggl(\frac{1}{n} \sum_{i=1}^n
x^{(i)}_{\setminus r} \bigl(x^{(i)}_{\setminus r}\bigr)^T\Biggr) \geq
D_{\max}\Biggr]\\
&&\qquad\leq
\mathbb{P} \Biggl[\Biggl\Vert\hspace*{-1.9pt}\Biggl\vert\Biggl(\frac{1}{n} \sum
_{i=1}^nx^{(i)}_{\setminus r} \bigl(x^{(i)}_{\setminus r}\bigr)^T \Biggr) -
{\mathbb{E}}_{{\theta^{*}}}[X_{\setminus r} X_{\setminus r}^T]
\Biggl\Vert\hspace*{-1.9pt}\Biggl\vert_{{2}} \geq
\delta\Biggr],
\end{eqnarray*}
which obeys the same upper bound (\ref{EqnUseful}) by following the
analogous argument.


\subsection[Proof of Lemma 6]{Proof of Lemma \protect\ref{LemConcentrateInco}}\label{sec52}

We begin by decomposing the sample matrix as the sum $ Q^{n}_{{S^{c}}
S} (Q^{n}_{SS})^{-1} = T_1 + T_2 + T_3 +
T_4$ where we define
%
%
\begin{subequation}
\begin{eqnarray}
T_1 & :=& Q^{*}_{{S^{c}}S} [ (Q^{n}_{S
S})^{-1} - (Q^{*}_{SS})^{-1} ] ,\\
T_2 & :=& [Q^{n}_{{S^{c}}S} - Q^{*}_{{S^{c}}S}
] (Q^{*}_{SS})^{-1},\\
T_3 & :=& [Q^{n}_{{S^{c}}S}
- Q^{*}_{{S^{c}}S} ] [ (Q^{n}_{SS})^{-1} -
(Q^{*}_{SS})^{-1} ], \\
T_4 & :=& Q^{*}_{{S^{c}}S} (Q^{*}_{SS})^{-1}.
\end{eqnarray}
\end{subequation}
The fourth term is easily controlled; indeed, we have
\[
\nnorm T_4 \nnorm_{{\infty}} = \nnorm Q^{*}_{{S^{c}}S}
(Q^{*}_{SS})^{-1} \nnorm_{{\infty}} \leq1 - \alpha
\]
by the
incoherence assumption \textup{(A2)}. If we can show that
$\nnorm T_i \nnorm_{{\infty}} \leq\frac{\alpha}{6}$ for the
remaining indices $i=1,2,3$, then by our four term decomposition and
the triangle inequality, the sample version satisfies the
bound (\ref{EqnSampleMutinc}), as claimed. We deal with these
remaining terms using the following lemmas:
\begin{lems}\label{lemtechbounds}
For any $\delta> 0$ and constants $K, K'$, the following
bounds hold:
%
%
\begin{subequation}
\label{EqnTechBound}
\begin{eqnarray}
\label{EqnTechBoundA}
&&\mathbb{P}[\nnorm Q^{n}_{{S^{c}}S} - Q^{*}_{{S^{c}} S} \nnorm
_{{\infty}} \geq\delta] \nonumber\\[-8pt]\\[-8pt]
&&\qquad \leq 2 \exp\biggl(
-K
\frac{n \delta^2}{d^2} + \log(d) + \log
(p-
d) \biggr); \nonumber\\
\label{EqnTechBoundB}
&&\mathbb{P}[\nnorm Q^{n}_{SS} - Q^{*}_{S S} \nnorm_{{\infty}} \geq
\delta] \nonumber\\[-8pt]\\[-8pt]
&&\qquad \leq 2 \exp\biggl(-K
\frac{n \delta^2}{d^2} + 2 \log(d) \biggr);\nonumber
\\
\label{EqnTechBoundC}
&&\mathbb{P}[\nnorm(Q^{n}_{SS})^{-1} - (Q^{*}_{S S})^{-1} \nnorm
_{{\infty}} \geq\delta] \nonumber\\[-8pt]\\[-8pt]
&&\qquad \leq 4 \exp\biggl(
-K\frac{n \delta^2}{d^3} + K' \log
(d)\biggr).\nonumber
\end{eqnarray}
\end{subequation}
\end{lems}

See Appendix \ref{AppTechBounds} for the proof of these
claims.

\subsubsection*{Control of first term}

Turning to the first term, we
re-factorize it as
\[
T_1 = Q^{*}_{{S^{c}}S} (Q^{*}_{SS})^{-1}
[Q^{n}_{SS} - Q^{*}_{SS} ]
(Q^{n}_{SS})^{-1}
\]
and then bound it (using the sub-multiplicative property
$\nnorm A B \nnorm_{{\infty}} \leq\break\nnorm A \nnorm_{{\infty}}
\nnorm B \nnorm_{{\infty}}$) as follows:
\begin{eqnarray*}
\nnorm T_1 \nnorm_{{\infty}} & \leq& \nnorm Q^{*}_{{S^{c}} S}
(Q^{*}_{SS})^{-1} \nnorm_{{\infty}}
\nnorm Q^{n}_{SS} - Q^{*}_{SS} \nnorm_{{\infty}}
\nnorm(Q^{n}_{SS})^{-1} \nnorm_{{\infty}} \\
& \leq& (1-\alpha) \nnorm Q^{n}_{SS} - Q^{*} _{S S} \nnorm
_{{\infty}} \bigl\{ \sqrt{d} \nnorm(Q^{n}_{S S})^{-1} \nnorm
_{{2}} \bigr\},
\end{eqnarray*}
where we have used the incoherence assumption \textup{(A2)}. Using the
bound (\ref{EqnEigConcMin}) from Lemma \ref{LemConcentrateEig} with
$\delta=
C_{\min}/2$, we have $\nnorm(Q^{n}_{SS})^{-1} \nnorm_{{2}} =
[\Lambda_{\min}(Q^{n}_{SS})]^{-1} \leq\frac{2}{C_{\min}}$ with
probability greater than $1-\exp(-Kn/ d^2+ 2
\log(d) )$. Next, applying the bound (\ref{EqnTechBoundB})
with $\delta= c/\sqrt{d}$, we conclude that with probability
greater than $1-2 \exp( -Knc^2/d^3 +
\log(d) )$, we have
\[
\nnorm Q^{n}_{SS} - Q^{*}_{S S} \nnorm_{{\infty}} \leq c/\sqrt{d}.
\]
By choosing the constant $c > 0$ sufficiently small, we are guaranteed
that
%
%
\begin{equation}
\mathbb{P}[ \nnorm T_1 \nnorm_{{\infty}} \geq\alpha/6] \leq
2\exp\biggl( -K\frac{nc^2}{d^3} + \log(d) \biggr).
\end{equation}

\subsubsection*{Control of second term}

To bound $T_2$, we first write
\begin{eqnarray*}
\nnorm T_2 \nnorm_{{\infty}} & \leq& \sqrt{d}
\nnorm(Q^{*}_{SS})^{-1} \nnorm_{{2}} \nnorm Q^{n} _{{S^{c}} S} -
Q^{*}_{{S^{c}}S} \nnorm_{{\infty}} \\
& \leq& \frac{\sqrt{d}}{C_{\min}} \nnorm Q^{n}_{{S^{c}} S} -
Q^{*}_{{S^{c}}S} \nnorm_{{\infty}}.
\end{eqnarray*}
We then apply bound (\ref{EqnTechBoundA}) with $\delta=
\frac{\alpha}{3} \frac{C_{\min}}{\sqrt{d}}$ to conclude that
%
%
\begin{equation}
\mathbb{P}[\nnorm T_2 \nnorm_{{\infty}} \geq\alpha/3] \leq2
\exp\biggl( -K\frac{n}{d^3} + \log(p-d) \biggr).
\end{equation}

\subsubsection*{Control of third term}

Finally, in order to bound the third term $T_3$, we apply the bounds
(\ref{EqnTechBoundA}) and (\ref{EqnTechBoundB}), both with $\delta=
\sqrt{\alpha/3}$, and use the fact that $\log(d) \leq\log(p- d)$ to
conclude that
%
%
\begin{equation}
\mathbb{P}[ \nnorm T_3 \nnorm_{{\infty}} \geq\alpha/3] \leq4
\exp\biggl( -K\frac{n}{d^3} + \log(p- d) \biggr).
\end{equation}

Putting together all of the pieces, we conclude that
\[
\mathbb{P}[ \nnorm Q^{n}_{{S^{c}}S} (Q^{n}_{S S})^{-1} \nnorm
_{{\infty}} \geq1 - \alpha/2] = {\mathcal{O}}
\biggl(\exp\biggl(-K\frac{n}{d^3} + \log(p) \biggr) \biggr)
\]
as claimed.

%
\section{Experimental results}
\label{SecExperiments}

We now describe experimental results that illustrate some consequences
of Theorem \ref{ThmMain}, for various types of graphs and scalings of
$(n, p, d)$. In all cases, we solved the
$\ell_1$-regularized logistic regression using special purpose
interior-point code developed by Koh, Kim and Boyd \cite{KohKimBoy07}.

\begin{figure}

\includegraphics{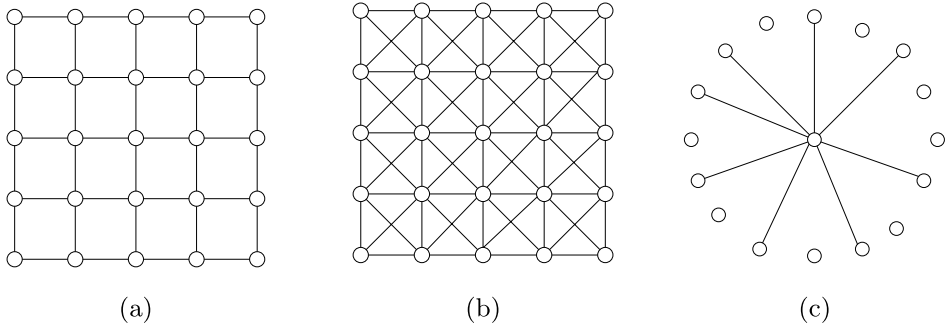}

\caption{Illustrations of different graph classes used in simulations.
\textup{(a)} Four-nearest neighbor grid ($d= 4$). \textup{(b)} Eight-nearest
neighbor grid ($d= 8$). \textup{(c)} Star-shaped graph
[$d=\Theta(p)$, or $d= \Theta(\log(p))$].}
\label{FigGraphs}
\end{figure}

We performed experiments for three different classes of graphs:
four-nearest neighbor lattices, (b) eight-nearest neighbor lattices
and (c) star-shaped graphs as illustrated in Figure \ref{FigGraphs}.
Given a distribution $\mathbb{P}_{{{\theta^{*}}}}$ of the Ising
form (\ref{EqnIsing}), we generated random data sets $\{x^{(1)},
\ldots, x^{(n)} \}$ by Gibbs sampling for the lattice models,
and by exact sampling for the star graph. For a given graph class and
edge strength $\omega> 0$, we examined the performance of models with
\textit{mixed couplings} meaning that $\theta^*_{st} = \pm\omega$ with
equal probability or with \textit{positive couplings} meaning that
$\theta^*_{st} = \omega$ for all edges $(s,t)$. In all cases, we set
the regularization parameter $\lambda_n$ as a constant factor
of $\sqrt{\frac{\log p}{n}}$
as suggested by Theorem \ref{ThmMain}. For any given graph and
coupling type, we performed simulations for sample sizes $n$
scaling as $n= 10 \beta d\log(p)$ where the
control parameter $\beta$ ranged from $0.1$ to upwards of $2$,
depending on the graph type.

%
%
%

%
%
\begin{figure}[b]

\includegraphics{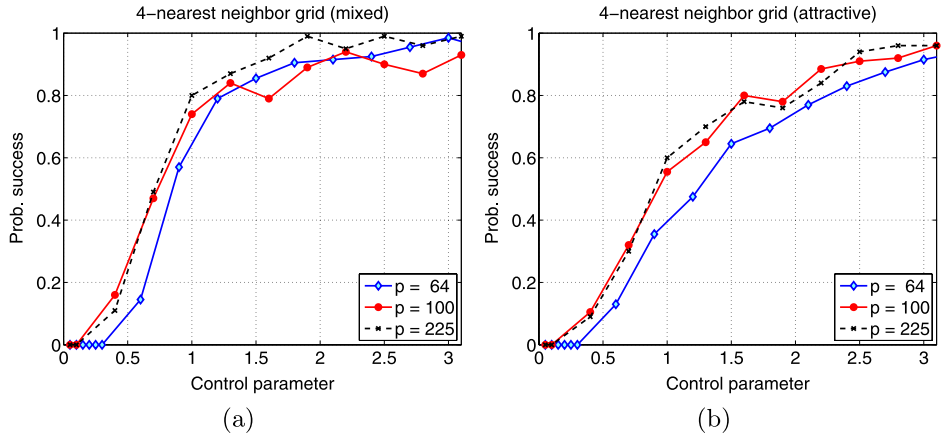}

\caption{Plots of success probability $\mathbb{P}[\widehat{\mathcal
{N}}_\pm(r) =
\mathcal{N}(r), \forall r]$ versus the control parameter $\beta
(n,
p, d) = n/[10 d\log(p)]$ for Ising models
on 2-D grids with four nearest-neighbor interactions ($d= 4$).
\textup{(a)} Randomly chosen mixed sign couplings $\theta^*_{st} = \pm0.50$.
\textup{(b)} All positive couplings $\theta^*_{st} = 0.50$.}
\label{FigIsingFour}
\end{figure}

\begin{figure}

\includegraphics{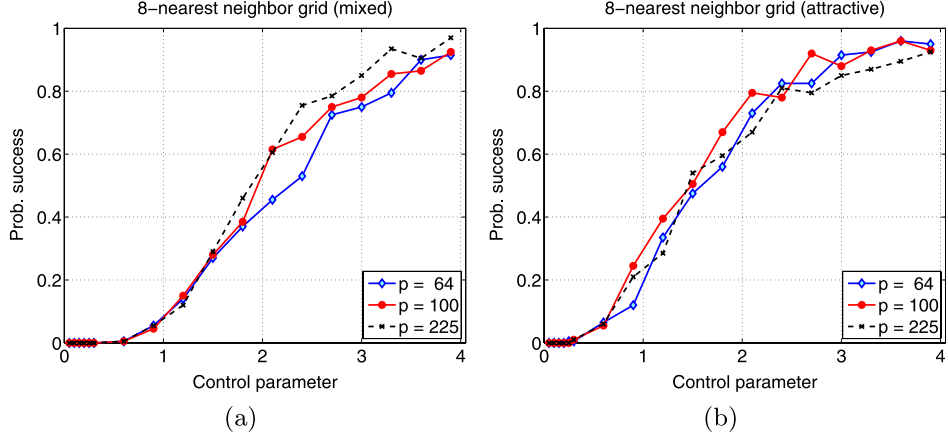}

\caption{Plots of success probability $\mathbb{P}[\widehat{\mathcal
{N}}_\pm(r) =
\mathcal{N}(r), \forall r]$ versus the control parameter $\beta
(n,p, d) = n/[10 d\log(p)]$ for Ising models
on 2-D grids with eight nearest-neighbor interactions ($d= 8$).
\textup{(a)} Randomly chosen mixed sign couplings $\theta^*_{st} = \pm0.25$.
\textup{(b)} All positive couplings $\theta^*_{st} = 0.25$.}
\label{FigIsingEight}
\end{figure}

Figure \ref{FigIsingFour} shows results for the $4$-nearest-neighbor
grid model, illustrated in Figure \ref{FigGraphs}(a) for three different
graph sizes $p\in\{64, 100, 225\}$ with mixed couplings [panel
(a)] and attractive couplings [panel (b)]. Each curve corresponds to
a given problem size, and corresponds to the success probability
versus the control parameter $\beta$. Each point corresponds to
the average of $N = 200$ trials. Notice how, despite the very
different regimes of $(n, p)$ that underlie each curve, the
different curves all line up with one another quite well. This fact
shows that for a fixed degree graph (in this case $\deg= 4$), the
ratio $n/\log(p)$ controls the success/failure of our model
selection procedure which is consistent with the prediction of
Theorem \ref{ThmMain}.
Figure \ref{FigIsingEight} shows analogous results for the
$8$-nearest-neighbor lattice model ($d= 8$), for the same range
of problem size $p\in\{64, 100, 225\}$ and for both mixed
and attractive couplings. Notice how once again the curves for
different problem sizes are all well aligned which is consistent with the
prediction of Theorem \ref{ThmMain}.

%

%
\begin{figure}[b]

\includegraphics{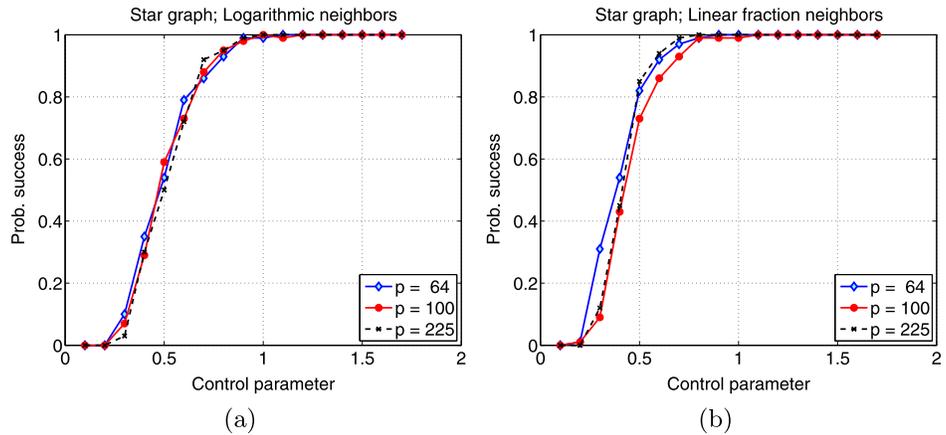}
 \caption{Plots of success probability
$\mathbb{P}[\widehat{\mathcal {N}}_\pm(r) = \mathcal{N}(r), \forall r]$
versus the control parameter $\beta (n, p, d) = n/[10 d\log(p)]$ for
star-shaped graphs for attractive couplings with \textup{(a)}
logarithmic growth in degrees, \textup{(b)} linear growth in degrees.}
\label{FigStar}
\end{figure}

For our next set of experiments, we investigate the performance
of our method for a class of graphs with unbounded maximum degree
$d$. In particular, we construct star-shaped graphs with
$p$ vertices by designating one node as the hub and connecting
it to $d< (p-1)$ of its neighbors. For linear
sparsity, we choose $d= \lceil0.1 p\rceil$, whereas for
logarithmic sparsity we choose $d= \lceil\log(p) \rceil$.
We again study a triple of graph sizes $p\in\{64, 100, 225 \}$,
and Figure \ref{FigStar} shows the resulting curves of success
probability versus control parameter $\beta= n/[10 d
\log(p)]$. Panels (a) and (b) correspond, respectively, to the
cases of logarithmic and linear degrees. As with the bounded degree
models in Figure \ref{FigIsingFour} and \ref{FigIsingEight}, these curves
align with one another showing a transition from failure to success
with probability one.

Although the purpose of our experiments is mainly to illustrate the
consequences of Theorem \ref{ThmMain}, we also include a comparison of
our nodewise $\ell_1$-penalized logistic regression-based method to two
other graph estimation procedures. For the comparison, we use a
star-shaped graph as in the previous plot, with one node designated as
the hub connected to $d= \lceil0.1 p\rceil$ of its
neighbors. It should be noted that among all graphs with a fixed
total number of edges, this class of graphs is among the most
difficult for our method to estimate. Indeed, the sufficient
conditions of Theorem \ref{ThmMain} scale logarithmically in the graph
size $p$ but polynomially in the maximum degree $d$;
consequently, for a fixed total number of edges, our method requires
the most samples when all the edges are connected to the same node, as in
a star-shaped graph.

For comparative purposes, we also illustrate the performance of the PC
algorithm of Spirtes, Glymour and Scheines \cite{spirtes00} as well as
the maximum
weight tree method of Chow and Liu \cite{chowliu68}. Since the star
graph is a tree (cycle-free), both of these methods are applicable in
this case. The PC algorithm is targeted to learning (equivalence
classes of) directed acyclic graphs, and consists of two stages. In
the first stage it starts from a completely connected undirected
graph, and iteratively removes edges based on conditional independence
tests so that at the end of this stage it is left with an undirected
graph which is called a skeleton. In the second stage, it partially
directs some of the edges in the skeleton so as to obtain a completed
partially directed acyclic graph which corresponds to an equivalence
class of directed acyclic graphs. As pointed out by Kalisch and
B\"{u}hlmann \cite{kalisch07}, for high-dimensional problems, the
output of the first stage, which is the undirected skeleton graph,
could provide a useful characterization of the dependencies in the
data. Following this suggestion, we use the skeleton graph determined
by the first stage of the PC algorithm as an estimate of the graph
structure. We use the \texttt{pcalg} R-package \cite{kalisch07} as
an implementation of the PC algorithm which uses partial
correlations to test conditional independencies.

%
\begin{figure}

\includegraphics{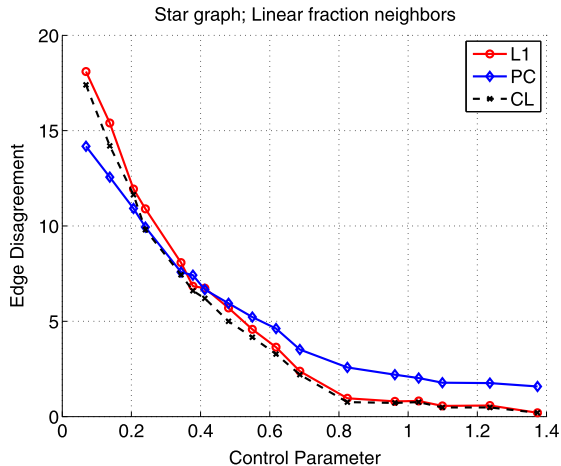}

\caption{Plots of edge disagreement $\mathbb{E}[ |\{ (s,t)
\mid
\widehat{E}_{st} \neq E^*_{st} \} |]$ versus the control parameter
$\beta(n, p, d) = n/[10 d
\log(p)]$ for star-shaped graphs where the hub node has degree
$d= \Theta(p)$. The results here are shown for attractive
couplings with $\theta^*_{st} = 0.25$ for all edges $(s,t)$ belonging
to the edge set. The $\ell_1$-penalized logistic regression method
(L1), the PC method (PC) and the maximum weight forest method of Chow
and Liu (CL) are compared for $p= 64$.}
\label{FigEdgeDis}
\end{figure}

The Chow--Liu algorithm \cite{chowliu68} is a method for exact
maximum likelihood structure selection which is applicable to the case of
trees. More specifically, it chooses, from among all trees with a
specified number of edges, the tree that minimizes the
Kullback--Leibler divergence to the empirical distribution defined by
the samples. From an implementational point of view, it starts with a
completely connected weighted graph with edge weights equal to the
empirical mutual information between the incident node variables of
the edge and then computes its maximum weight spanning tree. Since
our underlying model is a star-shaped graph with fewer than
$(p-1)$ edges, a spanning tree would necessarily include false
positives. We thus estimate the maximum weight forest with $d$
edges instead where we supplied the number of edges $d$ in the
true graph to the algorithm.

Figure \ref{FigEdgeDis} plots, for the three methods, the total number
of edge disagreements between the estimated graphs and the true graph
versus the control parameter $\beta= n/[10 d
\log(p)]$. Even though this class of graphs is especially
challenging for a neighborhood-based method, the $\ell_1$-penalized
logistic-regression based method is competitive with the Chow--Liu
algorithm, and except at very small sample sizes, it performs better
than the PC algorithm for this problem.

\section{Extensions to general discrete Markov random fields}
\label{SecGenMRF}

Our method and analysis thus far has been specialized to the case of
the binary pairwise Markov random fields. In this section, we briefly
outline the extension to the case of general discrete pairwise Markov
random fields. (Recall that for discrete Markov random fields, there
is no loss of generality in assuming only pairwise interactions since
by introducing auxiliary variables, higher-order interactions can be
reformulated in a pairwise manner \cite{WaiJor03Monster}.) Let $X =
(X_1,\ldots,X_p)$ be a random vector, each variable $X_i$ taking
values in a set $\mathcal{X}$ of cardinality $m$, say $\mathcal{X} =
\{1,2, \ldots, {m}\}$. Let $G= (V,E)$ denote a
graph with $p$ nodes corresponding to the $p$ variables
$\{X_1,\ldots,X_p\}$, and let $\{\phi_s\dvtx\mathcal{X} \rightarrow
{\mathbb{R}}, s \in V\}$ and $\{\phi_{st}\dvtx\mathcal{X} \times
\mathcal{X} \rightarrow\mathbb{R}, (s,t) \in E\}$, respectively,
denote collections of potential functions associated with the nodes
and edges of the graph. These functions can be used to define a
pairwise Markov random field over $(X_1, \ldots, X_p)$, with
density
%
%
\begin{equation}
\label{EqnGenMRF}
\mathbb{P}(x) \propto\exp\biggl\{\sum_{s \in V}
\phi_s(x_s) + \sum_{(s,t) \in E} \phi_{st}(x_s,x_t) \biggr\}.
\end{equation}

Since $\mathcal{X}$ is discrete, each potential function $\phi_{st}$
can be parameterized as linear combinations of $\{0,1\}$-valued
indicator functions. In particular, for each $s \in V$ and $j
\in\{1, \ldots, {m}-1\}$, we define
\[
\mathbb{I}[x_{s} = j] = \cases{1, &\quad if $x_s =j$,\cr
0, &\quad otherwise.}
\]
Note we omit an indicator for $x_s = {m}$ from the list, since
it is redundant given the indicators for $j = 1, \ldots, {m}-1$.
In a similar fashion, we define the pairwise indicator functions
$\mathbb{I}[x_s =j, x_t =k]$, for $(j,k) \in\{1, 2, \ldots, {m}-1\}^2$.

Any set of potential functions can then be written as
\[
\phi_{s}(x_s) = \sum_{j \in\{1, \ldots, {m}-1 \}}
\theta^*_{s;j} \mathbb{I}[x_s = j] \qquad\mbox{for $s \in
V$},
\]
and
\[
\phi_{st}(x_s,x_t) = \sum_{(j,k) \in\{1, \ldots, {m}-1
\}^2} \theta^*_{st;jk} \mathbb{I}[x_s = j, x_t = k]\qquad
\mbox{for $(s,t) \in E$.}
\]
Overall, the Markov random field can be parameterized in terms of the
vector $\theta^*_s \in{\mathbb{R}}^{{m}-1}$ for each $s \in V$,
and the vector $\theta^*_{st} \in{\mathbb{R}}^{({m}-1)^2}$ associated
with each edge. In discussing graphical model selection, it is
convenient to associate a vector $\theta^*_{uv} \in
{\mathbb{R}}^{({m}-1)^2}$ to \textit{every} pair of distinct vertices
$(u,v)$ with the understanding that $\theta^*_{uv} = 0$ if $(u,v)
\notin E$.

With this set-up, we now describe a graph selection procedure that is
the natural generalization of our procedure for the Ising model. As
before we focus on recovering for each vertex $r\in V$ its
neighborhood set and then combine the neighborhood sets across
vertices to form the graph estimate.

For a binary Markov random field (\ref{EqnIsing}), there is a unique
parameter $\theta^*_{rt}$ associated with each edge $(r, t)
\in E$. For ${m}$-ary models, in contrast, there is a
vector $\theta^*_{rt} \in{\mathbb{R}}^{({m}-1)^2}$ of parameters
associated with any edge $(r,t)$. In order to describe a
recovery procedure for the edges, let us define a matrix
$\Theta^*_{\setminus r}\in{\mathbb{R}}^{({m}-1)^2 \times
(p-1)}$ where column $u$ is given by the vector $\theta^*_{r
u}$. Note that unless vertex $r$ is connected to all of its
neighbors, many of the matrix columns are zero. In particular, the
problem of neighborhood estimation for vertex $r$ corresponds to
estimating the \textit{column support} of the matrix
$\Theta^*_{\setminus r}$---that is,
\[
\mathcal{N}(r) = \bigl\{ u \in V\setminus\{r\}
\mid\|\theta^*_{ru}\|_2 \neq0 \bigr\}.
\]

In order to estimate this column support, we consider the conditional
distribution of $X_r$ given the other variables $X_{\setminus
\{r\}} = \{X_{t} | t \in V \setminus\{r\}\}$. For a
binary model, this distribution is of the logistic form while for a
general pairwise MRF, it takes the form
%
%
\begin{equation}\qquad
\label{EqnMultiClassLogistic}
\mathbb{P}_\Theta[X_r= j \mid X_{\setminus r} =
x_{\setminus r} ] = \frac{\exp(\theta_{r;j}^* +
\sum_{t \in V\setminus\{r\}} \sum_{k} \theta^*_{rt;jk}
\mathbb{I}[x_t = k] )} {\sum_{\ell}\exp(\theta_{r;\ell
}^* +
\sum_{t \in V\setminus\{r\}} \sum_{k} \theta^*_{r
t;\ell
k} \mathbb{I}[x_t = k] )}.
\end{equation}
Thus, $X_r$ can be viewed as the response variable in a multiclass
logistic regression in which the indicator functions associated with
the other variables,
\[
\bigl\{\mathbb{I}[x_t = k] , t \in V \setminus\{r\}, k\in\{1, 2,
\ldots,
{m}-1\} \bigr\},
\]
play the role of the covariates.

Accordingly, one method of recovering the row support of
$\Theta^*_{\setminus r}$ is by performing multiclass logistic
regression of $X_r$ on the rest of the variables $X_{\setminus
r}$ using a block $\ell_2/\ell_1$ penalty of the form
\[
\nnorm\Theta_{\setminus r} \nnorm_{{2,1}} := \sum_{u \in
V\setminus\{r\}} \|\theta_{ru}\|_2.
\]
More specifically, let $\mathfrak{X}_1^n= \{x^{(1)},\ldots,
x^{(n)}\}$ denote an i.i.d. set of $n$ samples, drawn
from the discrete MRF (\ref{EqnGenMRF}). In order to estimate the
neighborhood of node $r$, we solve the following convex program:
%
%
\begin{equation}
\label{EqnGroup}
\min_{\Theta_{\setminus r} \in{\mathbb{R}}^{({m}-1)^2 \times
(p-1)}} \{ \ell(\Theta_{\setminus r}; \mathfrak{X}_1^n) +
\lambda_n\nnorm\Theta_{\setminus r} \nnorm_{{2,1}}\},
\end{equation}
where $\ell(\Theta_{\setminus r}; \mathfrak{X}_1^n) :=
\frac{1}{n} \sum_{i=1}^n\log\mathbb{P}_{\Theta}
[x^{(i)}_{r} \mid x^{(i)}_{\setminus r} ]$
is the rescaled multiclass logistic likelihood defined by the
conditional distribution (\ref{EqnMultiClassLogistic}), and
$\lambda_n> 0$ is a regularization parameter.

The convex program (\ref{EqnGroup}) is the multiclass logistic analog
of the group Lasso, a~type of relaxation that has been studied in
previous and on-going work on linear and logistic regression
(e.g., \cite{Kim05,meier,OboWaiJor08,YuaLi06}). It should be possible
to extend our analysis from the preceding sections so as to obtain
similar high-dimensional consistency rates for this multiclass
setting; the main difference is the slightly different
sub-differential associated with the block $\ell_2/\ell_1$ norm. See
Obozinski, Wainwright and Jordan \cite{OboWaiJor08} for some related
work on support
recovery using $\ell_2/\ell_1$ block-regularization for multivariate
linear regression.

\section{Conclusion}
\label{SecDiscussion}

We have shown that a technique based on $\ell_1$-regularized logistic
regression can be used to perform consistent model selection in binary
Ising graphical models, with polynomial computational complexity and
sample complexity logarithmic in the graph size. Our analysis applies
to the high-dimensional setting, in which both the number of nodes
$p$ and maximum neighborhood sizes $d$ are allowed to grow
as a function of the number of observations $n$. Simulation
results show the accuracy of these theoretical predictions. For
bounded degree graphs, our results show that the structure can be
recovered with high probability once $n/\log(p)$ is
sufficiently large. Up to constant factors, this result matches known
information-theoretic lower bounds \cite{SanWai08}. Overall, our
experimental results are consistent with the conjecture that logistic
regression procedure fails with high probability for sample sizes
$n$ that are smaller than ${\mathcal{O}}(d\log p)$. It
would be interesting to prove such a converse result, to parallel the
known upper and lower thresholds for success/failure of
$\ell_1$-regularized linear regression, or the Lasso
(see~\cite{Wainwright06new}).

As discussed in Section \ref{SecGenMRF}, although the current analysis
is applied to binary Markov random fields, the methods of this paper
can be extended to general discrete graphical models with a higher
number of states using a multinomial likelihood and some form of
block regularization. It should also be possible and would be
interesting to obtain high-dimensional rates in this setting. A final
interesting direction for future work is the case of samples drawn in
a non-i.i.d. manner from some unknown Markov random field; we suspect
that similar results would hold for weakly dependent sampling schemes.


\begin{appendix}
\section{Proof of uniqueness lemma}
\label{AppLemUnique}

In this appendix, we prove Lemma \ref{LemUnique}. By Lagrangian
duality, the penalized problem (\ref{EqnRegLikeTwo}) can be written as
an equivalent constrained optimization problem over the ball
$\|\theta\|_1 \leq C(\lambda_n)$, for some constant
$C(\lambda_n) < +\infty$. Since the Lagrange multiplier
associated with this constraint---namely, $\lambda_n$---is
strictly positive, the constraint is active at any optimal solution
so that $\|\theta\|_1$ is constant across all optimal solutions.

By the definition of the $\ell_1$-subdifferential, the subgradient
vector $\widehat{z}$ can be expressed as a convex combination
of sign vectors of the form
%
%
\begin{equation}
\widehat{z} = \sum_{v\in\{-1,+1\}^{p-1}} \alpha
_v v,
\end{equation}
where the weights $\alpha_v$ are nonnegative and sum to one.
In fact, these weights correspond to an optimal vector of Lagrange
multipliers for an alternative formulation of the problem in which
$\alpha_v$ is the Lagrange multiplier for the constraint
$\langle v , \theta\rangle\leq C(\lambda_n)$. From standard
Lagrangian theory \cite{Bertsekasnonlin}, it follows that any other
optimal primal solution $\widetilde{\theta}$ must minimize the associated
Lagrangian---or equivalently, satisfy
(\ref{EqnZeroGrad})---and moreover must satisfy the
complementary slackness conditions $\alpha_v
\{\langle v , \widetilde{\theta} \rangle- C \} = 0$ for all\vspace*{1pt} sign
vectors $v$. But these conditions imply that
$\langle\widehat{z} , \widetilde{\theta} \rangle= C = \|
\widetilde{\theta}\|_1$ which
cannot occur if $\widetilde{\theta}_j \neq0$ for some index $j$ for which
$|\widehat{z}_j| < 1$. We thus conclude that $\widetilde{\theta}_{{{S^{c}}}}
= 0$ for all optimal primal solutions.

Finally, given that all optimal solutions satisfy $\theta_{{{S^{c}}}}
= 0$, we may consider the restricted optimization problem subject to
this set of constraints. If the principal submatrix of the Hessian is
positive definite, then this sub-problem is strictly convex so that
the optimal solution must be unique.


\section{Proofs for technical lemmas}
\label{AppTechnical}

In this section, we provide proofs of
Lemmas \ref{LemBernOne}, \ref{Leml2cons} and \ref{LemTayRem},
previously stated in Section \ref{SecFixedDesign}.

\subsection{\texorpdfstring{Proof of Lemma \protect\ref{LemBernOne}}{Proof of Lemma 2}}
\label{AppBernOne}

Note that any entry of $W^n$ has the form $W^n_u =
\frac{1}{n} \sum_{i=1}^nZ^{(i)}_u$ where for $i=1, 2,
\ldots, n$, the variables
\[
Z^{(i)}_u := x^{(i)}_{\setminus r} \bigl\{x^{(i)}_r
- \mathbb{P}_{{{\theta^{*}}}}\bigl[x_r= 1 \mid x^{(i)}_{\setminus r
}\bigr] +
\mathbb{P}_{{{\theta^{*}}}}\bigl[x_r= -1 \mid x^{(i)}_{\setminus r}\bigr]
\bigr\}
\]
are zero-mean under $\mathbb{P}_{{{\theta^{*}}}}$, i.i.d. and bounded
($|Z^{(i)}_u| \leq2$). Therefore, by the Azuma-Hoeffding
inequality \cite{Hoeffding63}, we have, for any $\delta> 0$, $\mathbb{P}
[|W^n_u| > \delta] \leq2 \exp( -
\frac{n\delta^2}{8} )$. Setting $\delta= \frac
{\alpha
\lambda_n}{4 (2-\alpha)}$, we obtain
\[
\mathbb{P} \biggl[\frac{2-\alpha}{\lambda_n} |W^n_u| >
\frac{\alpha}{4} \biggr] \leq2 \exp\biggl( - \frac{\alpha^2
\lambda_n^2}{128 (2-\alpha)^2} n \biggr).
\]
Finally, applying a union bound over the indices $u$ of $W^n$
yields
\[
\mathbb{P} \biggl[\frac{2-\alpha}{\lambda_n} \|W^n\|
_\infty>
\frac{\alpha}{4} \biggr] \leq2 \exp\biggl( - \frac{\alpha^2
\lambda_n^2}{128 (2-\alpha)^2} n+ \log(p) \biggr)
\]
as claimed.
%

\subsection{\texorpdfstring{Proof of Lemma \protect\ref{Leml2cons}}{Proof of Lemma 3}}
\label{Appl2cons}
Following a method used in a different context by Rothman et
al. \cite{Rothman08}, we define the function $G\dvtx{\mathbb{R}}^{d}
\rightarrow{\mathbb{R}}$ by
%
%
\begin{equation}
G(u_S) :=\ell(\theta^*_{S} + u_S;
\mathfrak{X}_1^n) - \ell(\theta^*_{S}; \mathfrak{X}_1^n) +
\lambda_n (\|\theta^*_{S} + u_S\| -
\|\theta^*_{S}\| ).
\end{equation}
It can be seen from (\ref{thetaSCon}) that $\widehat{u} =
\widehat{\theta}_{S} - \theta^*_{S}$ minimizes $G$. Moreover, $G(0)
= 0$ by
construction; therefore, we must have $G(\widehat{u}) \leq0$.
Note also that $G$ is convex. Suppose that we show that for some
radius $B > 0$, and for $u\in{\mathbb{R}}^d$ with $\|u\|_2 =
B$, we have $G(u) > 0$. We then claim that $\|\widehat{u}\|_2
\leq B$. Indeed, if $\widehat{u}$ lay outside the ball of radius
$B$, then the convex combination $t \widehat{u} + (1-t) (0)$ would
lie on the boundary of the ball, for an appropriately chosen $t \in
(0,1)$. By convexity,
\[
G \bigl(t \widehat{u} + (1-t) (0) \bigr)
\leq tG(\widehat{u}) + (1-t) G(0) \leq0,
\]
contradicting the assumed strict positivity of $G$ on the boundary.

It thus suffices to establish strict positivity of $G$ on the boundary
of the ball with radius $B = M\lambda_n
\sqrt{d}$ where $M> 0$ is a parameter to be chosen later
in the proof. Let $u\in{\mathbb{R}}^d$ be an arbitrary vector
with $\|u\|_2 = B$. Recalling the notation $W= \nabla
\ell({{\theta^{*}}}; \mathfrak{X}_1^n)$, by a Taylor series
expansion of the log likelihood component of~$G$, we have
%
%
\begin{equation}
\label{EqnGtaylor}
G(u) = W_{S}^T u+ u^T [ \nabla^2
\ell(\theta^*_{S} + \alpha u) ] u+
\lambda_n (\|\theta^*_{S} + u_{S}\| -
\|\theta^*_{S}\| )
\end{equation}
for some $\alpha\in[0,1]$.
For the first term, we have the bound
%
%
\begin{equation}
\label{EqnBoundOne}
|W_S^T u| \leq\|W_S\|_\infty\|u
\|_1
\leq\|W_S\|_\infty\sqrt{d} \|u\|_2
\leq\bigl(\lambda_n\sqrt{d} \bigr)^2
\frac{M}{4},
\end{equation}
since $\|W_S\|_\infty\leq\frac{\lambda_n}{4}$ by
assumption.

Applying the triangle inequality to the last term in the
expansion (\ref{EqnGtaylor}) yields
\[
\lambda_n\|\theta^*_{S} + u_{S}\|_{1} -
\|\theta^*_{S}\|_{1} \geq- \lambda_n\|u_S\|_1.
\]
Since $\|u_S\|_1 \leq\sqrt{d} \|u_S\|_2$,
we have
%
%
\begin{equation}
\label{EqnBoundTwo}
\lambda_n\|\theta^*_{S} + u_{S}\|_{1} -
\|\theta^*_{S}\|_{1} \geq-\lambda_n\sqrt{d}
\|u_S\|_2 = -M \bigl(\sqrt{d} \lambda
_n \bigr)^2.
\end{equation}

Finally, turning to the middle Hessian term, we have
\begin{eqnarray*}
q^* &:=& \Lambda_{\min}\bigl(\nabla^2\ell(\theta^*_{S} + \alpha u;
\mathfrak{X}_1^n)\bigr)\\
&\geq&\min_{\alpha\in[0,1]} \Lambda_{\min}\bigl(\nabla^2\ell
(\theta^*
_{S} +
\alpha u_{S}; \mathfrak{X}_1^n)\bigr)\\
&=& \min_{\alpha\in[0,1]}
\Lambda_{\min} \Biggl[\frac{1}{n}\sum_{i=1}^{n}\eta\bigl(x^{(i)}; \theta
^*_S + \alpha u_S\bigr) x^{(i)}_{S} \bigl(x^{(i)}_{S}\bigr)^T \Biggr].
\end{eqnarray*}
By a Taylor series expansion of $\eta(x^{(i)}; \cdot)$, we have
\begin{eqnarray*}
q^* &\geq&
\Lambda_{\min} \Biggl[\frac{1}{n}\sum_{i=1}^{n}\eta\bigl(x^{(i)}; \theta
^* _S\bigr)x^{(i)}_{S}
\bigl(x^{(i)}_{S}\bigr)^T \Biggr]\\
&&{} - \max_{\alpha\in[0,1]}
\Biggl\Vert\hspace*{-1.9pt}\Biggl\vert\frac{1}{n}\sum_{i=1}^{n}\eta'\bigl(x^{(i)}; \theta^*_S+ \alpha
u_S\bigr)\bigl(u_{S}^{T}x^{(i)}_{S}\bigr)x^{(i)}_{S} \bigl(x^{(i)}_{S}\bigr)^T \Biggl\Vert\hspace*{-1.9pt}\Biggl\vert_{{2}} \\
&=& \Lambda_{\min}(Q^{*}_{SS}) - \max_{\alpha\in[0,1]}
\Biggl\Vert\hspace*{-1.9pt}\Biggl\vert\frac{1}{n}\sum_{i=1}^{n}\eta'\bigl(x^{(i)}; \theta^*_S+ \alpha
u_S\bigr) \bigl(\bigl\langle u_{S} , x^{(i)}_{S} \bigr\rangle\bigr)
x^{(i)}_{S} \bigl(x^{(i)}_{S}\bigr)^T \Biggl\Vert\hspace*{-1.9pt}\Biggl\vert_{{2}} \\
&\geq&
C_{\min}- \max_{\alpha\in[0,1]}
\Biggl\Vert\hspace*{-1.9pt}\Biggl\vert\mathop{\underbrace{\frac{1}{n}\sum_{i=1}^{n}\eta'\bigl(x^{(i)}; \theta
^*_S+ \alpha u_S\bigr) \bigl(\bigl\langle u_{S} , x^{(i)}_{S} \bigr\rangle\bigr)
x^{(i)}_{S} \bigl(x^{(i)}_{S}\bigr)^T} \Biggl\Vert\hspace*{-1.9pt}\Biggl\vert_{{2}} }_{A(\alpha)\phantom{00}}.
\end{eqnarray*}

It remains to control the spectral norm of the matrices
$A(\alpha)$, for $\alpha\in[0,1]$. For any fixed $\alpha\in
[0,1]$ and $y \in{\mathbb{R}}^d$ with $\|y \|_2 = 1$, we have
\begin{eqnarray*}
\langle y , A(\alpha) y \rangle& = & \frac{1}{n}\sum_{i=1}^{n}
\eta'\bigl(x^{(i)}; \theta^*_S+ \alpha u_S\bigr)
\bigl[\bigl\langle u_{S} , x^{(i)}_{S} \bigr\rangle\bigr] \bigl[\bigl\langle
x^{(i)} _S , y \bigr\rangle\bigr]^2 \\
& \leq& \frac{1}{n}\sum_{i=1}^{n} \bigl|\eta'\bigl(x^{(i)}; \theta^*_S+
\alpha u_S\bigr) \bigr| \bigl| \bigl\langle u_{S} , x^{(i)}_{S} \bigr\rangle
\bigr|\bigl[\bigl\langle x^{(i)}_S , y \bigr\rangle\bigr]^2.
\end{eqnarray*}
Now note that $|\eta'(x^{(i)}; \theta^*_S+ \alpha u_S) |
\leq
1$, and
\[
\bigl|\bigl\langle u_{S} , x^{(i)}_{S} \bigr\rangle\bigr| \leq\sqrt{d}
\|u_{S}\|_2 = M\lambda_nd.
\]
Moreover, we have
\[
\frac{1}{n} \sum_{i=1}^{n} \bigl(\bigl\langle x^{(i)}_{S} , y
\bigr\rangle
\bigr)^2
\leq\Biggl\Vert\hspace*{-1.9pt}\Biggl\vert\frac{1}{n} \sum_{i=1}^{n}x^{(i)}_{S}\bigl(x^{(i)}_{S}\bigr)^T
\Biggl\Vert\hspace*{-1.9pt}\Biggl\vert_{{2}} \leq D_{\max}
\]
by assumption.
Combining these pieces, we obtain
\[
{\max_{\alpha\in[0,1]}} \nnorm A(\alpha) \nnorm_{{2}}
\leq D_{\max}M\lambda_nd \leq C_{\min}/2,
\]
assuming that $\lambda_n\leq\frac{C_{\min}}{2 MD_{\max}
d}$. We verify this condition momentarily, after we have
specified the constant $M$.

Under this condition, we have shown that
%
%
\begin{equation}
\label{EqnBoundThree}
q^* :=\Lambda_{\min}\bigl(\nabla^2\ell(\theta^*_{S} + \alpha u;
\mathfrak{X}_1^n
)\bigr) \geq C_{\min}/2.
\end{equation}
Finally, combining the bounds (\ref{EqnBoundOne}),
(\ref{EqnBoundTwo}), and (\ref{EqnBoundThree}) in the
expression (\ref{EqnGtaylor}), we conclude that
\[
G(u_S) \geq\bigl(\lambda_n\sqrt{d
} \bigr)^2\biggl\{ -\frac{1}{4} M+ \frac{C_{\min}}{2} M^2 - M
\biggr\}.
\]
This expression is strictly positive for $M= 5/C_{\min}$. Consequently,
as long as
\[
\lambda_n \leq\frac{C_{\min}}{2 MD_{\max}d} =
\frac{C_{\min}^2}{10 D_{\max}d}
\]
as assumed in the statement of the lemma, we are guaranteed that
\[
\|u_S\|_2 \leq M\lambda_n\sqrt{d}
=\frac{5}{C_{\min}} \lambda_n\sqrt{d}
\]
as claimed.

\subsection{\texorpdfstring{Proof of Lemma \protect\ref{LemTayRem}}{Proof of Lemma 4}}
\label{AppTayRem}

We first show that the remainder term $R^n$ satisfies the
bound $\|R^n\|_\infty\leq D_{\max}\|\widehat{\theta}_S-
\theta^*_S\|_2^2$. Then the result of
Lemma \ref{Leml2cons}---namely, that $\|\widehat{\theta}_S-
\theta^*_S\|_2 \leq\frac{5}{C_{\min}} \lambda_n
\sqrt{d}$---can be used to conclude that
\[
\frac{\|R^n\|_\infty}{\lambda_n} \leq\frac{25
D_{\max}}{C_{\min}^2} \lambda_nd
\]
as claimed in Lemma \ref{LemTayRem}.

Focusing on element $R^n_j$ for some index $j \in\{1,
\ldots, p\}$, we have
\begin{eqnarray*}
R^n_j & = & \bigl[\nabla^2 \ell\bigl(\onebar{\theta}^{(j)}; x\bigr)
- \nabla^2
\ell({{\theta^{*}}}; x) \bigr]_j^T [\widehat{\theta}- \theta
^*] \\
& = & \frac{1}{n} \sum_{i=1}^n
\bigl[\eta\bigl(x^{(i)}; \onebar{\theta}^{(j)}\bigr) - \eta\bigl(x^{(i)};
{{\theta^{*}}}\bigr)
\bigr] \bigl[x^{(i)}
\bigl(x^{(i)}\bigr)^T \bigr]_j^T [\widehat{\theta}- \theta^*]
\end{eqnarray*}
for some point $\onebar{\theta}^{(j)} = t_j\widehat{\theta}+
(1-t_j){{\theta^{*}}}$.
Setting $g(t) = \frac{4
\exp(2t)}{[1+ \exp(2t)]^2}$, note that $\eta(x; \theta) =
g(
x_r\sum_{t \in V\setminus r} \theta_{rt} x_t )$.
By the chain rule and another application of the
mean value theorem, we then have
\begin{eqnarray*}
R^n_j & = & \frac{1}{n} \sum_{i=1}^n
g' \bigl(\twobar{\theta}^{(j)T} x^{(i)} \bigr) \bigl(x^{(i)}\bigr)^T \bigl[\onebar
{\theta}
^{(j)} - \theta^*\bigr] \bigl\{
x^{(i)}_j \bigl(x^{(i)}\bigr)^T [\widehat{\theta}- \theta^*] \bigr\} \\
& = & \frac{1}{n} \sum_{i=1}^n
\bigl\{g' \bigl(\twobar{\theta}^{(j)T} x^{(i)} \bigr) x^{(i)}_j\bigr\}
\bigl\{\bigl[
\onebar{\theta}^{(j)}-{{\theta^{*}}}\bigr]^T x^{(i)}\bigl(x^{(i)}\bigr)^T
[\widehat{\theta}- \theta^*] \bigr\},
\end{eqnarray*}
where $\twobar{\theta}^{(j)}$ is another point on the line joining
$\widehat{\theta}$ and ${{\theta^{*}}}$.
Setting $a_i :=\{g' (\twobar{\theta}^{(j)T}\times x^{(i)} )
x^{(i)}_j \}$ and $b_i
:=\{ [\onebar{\theta}^{(j)}-\theta^*]^T x^{(i)}(x^{(i)})^T [\widehat
{\theta}-
\theta^*] \}$, we have
\[
|R^n_j| = \frac{1}{n} \Biggl|\sum_{i=1}^na_i
b_i \Biggr| \leq\frac{1}{n} \|a\|_\infty\|b\|_1.
\]
A calculation shows that $\|a\|_\infty\leq1$, and
\begin{eqnarray*}
\frac{1}{n} \|b\|_1 & = & t_j [\widehat{\theta}- \theta^*]^T
\Biggl\{\frac{1}{n} \sum_{i=1}^n x^{(i)}\bigl(x^{(i)}\bigr)^T \Biggr\} [\widehat
{\theta}-
\theta^*] \\
& = & t_j [\widehat{\theta}_S- \theta^*_S]^T \Biggl\{\frac
{1}{n}
\sum_{i=1}^n x^{(i)}_S\bigl(x^{(i)}_S\bigr)^T \Biggr\} [\widehat{\theta}_S-
\theta^*_S] \\
& \leq& D_{\max}\|\widehat{\theta}_S- \theta^*_S\|_2^2,
\end{eqnarray*}
where the second line uses the fact that $\widehat{\theta}_{S^{c}}=
\theta^*_{S^{c}}= 0$. This concludes the proof.

\section{\texorpdfstring{Proof of Lemma \protect\lowercase{\ref{lemtechbounds}}}{Proof of Lemma 7}}
\label{AppTechBounds}

Recall from the discussion leading up to the bound (\ref{EqnAzuma})
that element $(j,k)$ of the matrix difference $Q^{n}- Q^{*}$,
denoted by $Z_{jk}$, satisfies a sharp tail bound. By definition of
the $\ell_\infty$-matrix norm, we have
\begin{eqnarray*}
\mathbb{P}[\nnorm Q^{n}_{{S^{c}}S} - Q^{*}_{{S^{c}} S} \nnorm
_{{\infty}} \geq\delta] & = & \mathbb{P}\biggl[{\max_{j \in
{S^{c}}} \sum_{k
\in S} }|Z_{jk}| \geq\delta\biggr] \\
& \leq& (p-d) \mathbb{P}\biggl[{\sum_{k\in S}} |Z_{jk}|
\geq
\delta\biggr],
\end{eqnarray*}
where the final inequality uses a union bound, and the fact that
$|{S^{c}}| \leq p- d$. Via another union bound over the
row elements
\begin{eqnarray*}
\mathbb{P}\biggl[{\sum_{k \in S}} |Z_{jk}| \geq\delta\biggr] &\le&\mathbb{P}
[\exists k \in S | |Z_{jk}| \geq\delta/d]\\
&\le& d \mathbb{P}[|Z_{jk}| \geq\delta/d];
\end{eqnarray*}
we then obtain
\[
\mathbb{P}[\nnorm Q^{n}_{{S^{c}}S} - Q^{*}_{{S^{c}} S} \nnorm
_{{\infty}} \geq\delta] \leq(p- d)
d \mathbb{P} [|Z_{jk}| \geq\delta/d ]
\]
from which the claim (\ref{EqnTechBoundA}) follows by setting
$\varepsilon= \delta/d$ in the Hoeffding bound~(\ref{EqnAzuma}).
The proof of bound (\ref{EqnTechBoundB}) is analogous with the
pre-factor $(p- d)$ replaced by $d$.

To prove the last claim (\ref{EqnTechBoundC}), we write
\begin{eqnarray*}
\nnorm(Q^{n}_{SS})^{-1} - (Q^{*}_{S S})^{-1} \nnorm_{{\infty}} & = &
\nnorm(Q^{*}_{S S})^{-1} [Q^{*}_{SS} - Q^{n}_{SS} ]
(Q^{n}_{SS})^{-1} \nnorm_{{\infty}} \\
& \leq& \sqrt{d} \nnorm(Q^{*}_{S S})^{-1} [Q^{*}_{SS} -
Q^{n}_{SS} ] (Q^{n}_{SS})^{-1} \nnorm_{{2}} \\
& \leq& \sqrt{d} \nnorm(Q^{*}_{SS})^{-1} \nnorm_{{2}}
\nnorm Q^{*}_{SS} - Q^{n}_{SS} \nnorm_{{2}}
\nnorm(Q^{n}_{SS})^{-1} \nnorm_{{2}} \\
& \leq& \frac{\sqrt{d}}{C_{\min}} \nnorm Q^{*}_{S S} - Q^{n}_{SS}
\nnorm_{{2}} \nnorm(Q^{n}_{SS})^{-1} \nnorm_{{2}}.
\end{eqnarray*}
From the proof of Lemma \ref{LemConcentrateEig}, in particular
equation (\ref{EqnUseful}), we have
\[
\mathbb{P} \biggl[\nnorm(Q^{n}_{SS})^{-1} \nnorm_{{2}} \geq
\frac{2}{C_{\min}} \biggr] \leq2 \exp\biggl(-K\frac{\delta^2
n}{ d^2} + B\log(d) \biggr)
\]
for a constant $B$. Moreover, from (\ref{EqnUseful}),
we have
\[
\mathbb{P}\bigl[\nnorm Q^{n}_{SS} -Q_{SS} \nnorm_{{2}}
\geq
\delta/\sqrt{d} \bigr] \leq2 \exp\biggl(-K\frac
{\delta^2
n}{d^3} + 2 \log(d) \biggr)
\]
so that the bound (\ref{EqnTechBoundC}) follows.
\end{appendix}

\printaddresses

\end{document}